\documentclass[11pt]{article}
\usepackage{geometry}                
\usepackage{graphicx}
\usepackage{amssymb}
\usepackage{amsmath}
\usepackage{epstopdf}
\usepackage[all]{xypic}

\newtheorem{lemma}{Lemma}[section] 
\newtheorem{propos}[lemma]{Proposition}
\newtheorem{example}[lemma]{Example}
\newtheorem{theorem}[lemma]{Theorem}
\newtheorem{cor}[lemma]{Corollary}

\newtheorem{defin}[lemma]{Definition}

\newcommand{\C}{\mathbb{C}}

\newcommand{\Hom}{\mathrm{Hom}}

\newcommand{\proof}{{\noindent {\bfseries  Proof:}\quad }}

\newcommand{\ev}{\mathrm{ev}}
\newcommand{\coev}{\mathrm{coev}}

\newcommand{\extd}{\mathrm{d}}

\newcommand{\tens}{\mathop{\otimes}}
\newcommand{\la}{{\triangleright}}

\newcommand{\id}{\mathrm{id}}

\newcommand{\<}{\langle}
\renewcommand{\>}{\rangle}

\newcommand{\pdol}{\overline{\partial}}

\newcommand{\SU}{\mathrm{SU}}

\title{Differential and holomorphic differential operators on noncommutative algebras}
\author{Edwin Beggs}
\date{}                                           

\begin{document}

\maketitle

\abstract{This paper deals with sheaves of differential operators on noncommutative algebras. The sheaves are defined by quotienting a the tensor algebra of vector fields (suitably deformed by a covariant derivative) to ensure zero curvature. As an example we can obtain enveloping algebra like relations for Hopf algebras with differential structures which are not bicovariant. Symbols of differential operators are defined, but not studied. These sheaves are shown to be in the center of as category of bimodules with flat bimodule covariant derivatives. Also holomorphic differential operators are considered, though without the quotient to ensure zero curvature. }

\section{Introduction}
The purpose of this paper is to describe differential operators on differential graded algebras.
The zero grade of the differential graded algebra, $\Omega^0A$ or just $A$, is a possibly noncommutative algebra over the field of complex numbers, standing for a collection of smooth functions on a hypothetical `noncommutative space', and $\Omega^nA$ stands for the $n$-forms. Modules for $A$ stand for sections of vector bundles on the noncommutative space.

As we have no local coordinate patches, every time a textbook on differential operators would
mention partial derivatives, we
 have to use the corresponding globally defined object, a covariant derivative. However we have a complication, there is almost never a simple two sided Leibniz rule for noncommutative covariant derivatives. Every time where would have to use both a right and left Leibniz rule, we have to use a modification, a map $\sigma$ to reverse order of modules and 1-forms.
This idea had its origins in \cite{DVMic,DVMass} and \cite{Mourad}, and was later
used in \cite{Madore,FioMad}. In \cite{BMHDSfinitegrps} it was shown that this idea allowed tensoring
of bimodules with connections.
  This, and the usual complications of keeping track of order, mean that the proofs are more laborious than in the classical case. To make life even more difficult, we do not get a well defined object when applying a covariant derivative to one side of the tensor product $E\tens_A F$, as the
 covariant derivative is not a module map.
 In fact, many proofs are inductive, and very tedious.

In  \cite{bbdiff} noncommutative differential operators were defined in terms of actions of tensor products of vector fields on left modules with covariant derivative. This action required multiple derivatives, for which a bimodule covariant derivative on 1-forms on the algebra was needed.
By using this action, an associative algebra structure (called $\mathcal{T}\mathrm{Vec}A_\bullet$) involving differentiation was put on the tensor products of vector fields, and this gives composition of differential operators. If instead of acting on left modules with covariant derivative, we act on the moniodal category ${}_A\mathcal{E}_A$ of bimodules with left bimodule covariant derivative, we see that $\mathcal{T}\mathrm{Vec}A_\bullet$ is in the center of ${}_A\mathcal{E}_A$, and that the associated order reversing map gives the action on tensor products. This uses a canonical covariant derivative on $\mathcal{T}\mathrm{Vec}A_\bullet$, defined just as in the classical case.
Section~\ref{vgcyuai1} gives preliminaries on noncommutative covariant derivatives, and
Section~\ref{jkhsvcagkvc} summarises the work in \cite{bbdiff}. 

One thing not treated in \cite{bbdiff} is the idea of commutation of partial derivatives, i.e.\
$\frac{\partial }{\partial x}\, \frac{\partial }{\partial y}=\frac{\partial }{\partial y}\, \frac{\partial }{\partial x}$ on functions on $\mathbb{R}^2$. One reason for this is that this equality is not true in the context there - in terms of covariant derivatives the difference between the two sides is the curvature. However if we act just on the functions (or a module with zero curvature), we can impose the equality in classical geometry, and this is normally done in defining differential operators, to obtain the algebra usually known as $\mathcal{D}A$.
I have used \cite{MasSabDmod} as a reference for $\mathcal{D}$-modules. 
 Another way of looking at this classically is to consider the vector fields to be the Lie algebra associated to the diffeomorphism group, in which case the equality above (or rather its generalisation to arbitrary vector fields) is similar to the relation in the  universal enveloping algebra.

In many places in the literature the phrase `sheaf of differential operators' appears. In \cite{bbsheaf} a sheaf is defined in noncommutative differential geometry as a module with zero curvature covariant derivative. We can take the curvature of the canonical covariant derivative on $\mathcal{T}\mathrm{Vec}A_\bullet$, and this gives the curvature as a differential operator, just as in the classical case. Of course, this curvature is likely not zero. However, if the 2-forms are a finitely generated projective module, it turns out that there are relations which can be imposed on $\mathcal{T}\mathrm{Vec}A_\bullet$ to force the curvature to vanish, and we define $\mathcal{D}A$ to be $\mathcal{T}\mathrm{Vec}A_\bullet$ quotiented by these relations. The action of $\mathcal{T}\mathrm{Vec}A_\bullet$ on a module with zero curvature covariant derivative restricts to an action of $\mathcal{D}A$. In the classical case, we recover the commutation of partial derivatives above.

Naturally, recovering the classical case is not enough, and for an interesting example we take a noncommutative algebra. We shall take the example of the left covariant 3D differential calculus on deformed $SU(2)$ in \cite{W1,W2}. We can restrict attention to the left invariant vector fields on deformed
$SU(2)$ - classically this would be the Lie algebra - and write the relations in $\mathcal{D}A$ explicitly.
We shall, by choice of a covariant derivative from \cite{BMriem}, write down 3 by 3 matrices giving the action of $\mathcal{D}A$ on the left invariant 1-forms. The relations in $\mathcal{D}A$ \textit{look as though they were} a $q$ deformation of the universal enveloping algebra of $su(2)$.
Now Woronowicz in \cite{W2} did give a construction of a deformed Lie algebra \textit{in the bicovariant case}, which does not include the 3D calculus. Also Majid in \cite{MajBrLie}
introduced the universal enveloping algebra of a braided Lie algebra.
In fact there is a circle sub Hopf algebra of $SU(2)_q$ for which there is a right coaction on the 3D calculus, and a Yetter-Drinfeld braiding, but that is not the point - the construction of $\mathcal{D}A$ is independent of any bicovariance or braiding assumptions. As yet, it is not obvious what the proper algebraic interpretation of this should be, but it is essentially substituting knowledge about the 2-forms for bicovariance.

In Section~\ref{bvhdskv} we show that $\mathcal{D}A$  is in the center of the category ${}_A\mathcal{F}_A$ of bimodules with flat bimodule covariant derivatives. 
This is analagous to the result in \cite{bbdiff} for $\mathcal{T}\mathrm{Vec}A_\bullet$ being in the center of ${}_A\mathcal{E}_A$. In Section~\ref{bcadhklskuc} we consider combining
differential operators with module endomorphisms. This is done for two reasons: Firstly it removes any dependence of the theory on the choice of covariant derivative $\square$. Secondly it includes cases like the classical Dirac operator, which requires endomorphisms of the spinor bundle. In Section~\ref{bcdhlsku} we define symbols of noncommutative differential operators, but apart from giving an example to demonstrate that the corresponding polynomials are no longer commutative, say very little about this.

It is natural to ask what  of the theory of differential operators can be extended to 
noncommutative complex differential geometry. 
We shall use the version of noncommutative complex differential geometry from in 
\cite{BegSmiComplex} and referenced in \cite{KhLavS}, which is based on the classical approach set out in \cite{GH}. This is based on a bimodule map $J:\Omega^1A\to \Omega^1A$ satisfying $J^2$ being minus the identity and an integrability condition. Constructions of complex structures on noncommutative spaces has been carried out in various places, including
\cite{HK04,HK06,HK07,buncomproj}. Many aspects of the geometry of noncommutative projective space are discussed in \cite{DADL,DAL}, and the quantum plane is investigated in \cite{bdr}. Cocycle deformations of complex structures are discussed in \cite{bramaj}.
 In Section~\ref{bcudisokocfgy} we see that the definition of integrability for forms in \cite{BegSmiComplex} can be expressed in terms of vector fields, in a direct analogue of the classical Newlander-Nirenberg integrability theorem \cite{NewNir}

The first practical problem to be addressed in Section~\ref{ioygfcityk}, in the absence of local coordinates, is `what is a holomorphic vector field'. The simple answer that it is a vector field which, when applied to any holomorphic function, yields another holomorphic function. However this definition is not as simple as it seems, and requires some assumptions on covariant derivatives to give a simple answer. This is not surprising, as in the absence of local holomorphic coordinates, much of the description of the global complex structure falls on the shoulders of a covariant derivative on $\Omega^1 A$. Conditions on covariant derivatives are further discussed in Section~\ref{bcdsuovuy}. 

In Section~\ref{vukjhcfxfth} we see that the associative algebra $\mathcal{T}\mathrm{Vec}^{*}A_\bullet$ respects the $J$ structure, and that we have subalgebras $\mathcal{T}\mathrm{Vec}^{*,0}A_\bullet$ and $\mathcal{T}\mathrm{Vec}^{0,*}A_\bullet$ constructed from the $\pm\mathrm{i}$ eigenvectors of $J$ on $\mathrm{Vec}A$, and describe the categories of modules and covariant derivatives that they act on. In Section~\ref{nbcilhadskuv} we consider higher order derivatives of holomorphic functions, and classify those which preserve holomorphicity by 
yet another covariant derivative
$\bar\partial_{\mathcal{HD}}:\mathcal{T}\mathrm{Vec}^{*,0}A_\bullet\to \Omega^{0,1}A\tens_A \mathcal{T}\mathrm{Vec}^{*,0}A_\bullet$. It is shown that such holomorphic elements of 
$\mathcal{T}\mathrm{Vec}^{*,0}A_\bullet$ are preserved by the $\bullet$ product.

This paper stops short of defining the complex differential operators analogously to $\mathcal{D}A$ by quotienting $\mathcal{T}\mathrm{Vec}^{*,0}A_\bullet$ by an ideal to force zero holomorphic curvature. It is likely that this can be done, though at the possible cost of introducing more conditions on the covariant derivative $\square$, and this paper is long enough!

\section{Preliminaries on covariant derivatives} \label{vgcyuai1}
Here we use the same notation as in \cite{bbdiff}. 
Let $A$ be a unital algebra over $\C$. Suppose that the algebra $A$ has a
differential structure $(\Omega A,\extd,\wedge)$ in the sense of a differential
graded exterior algebra $\Omega A=\oplus_{n\ge 0}\Omega^n A$ with derivative $\extd:\Omega^nA\to \Omega^{n+1}A$ and product $\wedge:\Omega^nA\tens_A\Omega^mA\to \Omega^{n+m}A$. 
We have $\Omega^0A=A$, $\extd^2=0$ and 
 the graded Leibniz rule $\extd(\xi\wedge\eta)=\extd\xi\wedge\eta+(-1)^n\xi\wedge\extd\eta$ for 
$\xi\in\Omega^nA$.  We suppose that 
$\Omega^1A$ generates the exterior algebra over $A$, and that $\Omega^1 A=A.\extd A$. Note that we do not assume graded commutativity.

For a smooth manifold, the local coordinate patches give 1-forms which are finitely generated projective as a module over the functions on the manifold. In noncommutative geometry we do not have the coordinate patches, but the finitely generated projective assumption remains a sensible one to make.
We suppose that $\Omega^1 A$ is finitely generated projective as a right $A$ module, and set the vector fields $\mathrm{Vec} A=\Hom_A(\Omega^1 A,A)$ (i.e.\ the right module maps). We denote the corresponding evaluation and coevaluation maps by
\begin{eqnarray*}
\ev:\mathrm{Vec} A\tens_A \Omega^1 A\to A\ ,\quad \coev:A\to \Omega^1 A\tens_A \mathrm{Vec} A\ .
\end{eqnarray*}
(The coevaluation map is essentially the dual basis.)
For multiple copies,  define, where we have $n$ copies of $\mathrm{Vec} A$ and $\Omega^1 A$,
 \begin{eqnarray*}
&& \mathrm{Vec}^{\tens 0}A\ =\ A\ ,\quad
\mathrm{Vec}^{\tens n}A\ =\ \mathrm{Vec} A\tens_A \mathrm{Vec} A\tens_A\dots \tens_A \mathrm{Vec} A\ ,\cr
&& \Omega^{\tens 0}A\ =\ A\ ,\quad
\Omega^{\tens n}A\ =\ \Omega^1 A\tens_A \Omega^1 A\tens_A\dots \tens_A \Omega^1 A\ .
\end{eqnarray*}
Define the $n$-fold evaluation map $\ev^{\<n\>}:\mathrm{Vec}^{\tens n}A\tens_A \Omega^{\tens n}A\to A$
recursively by
\begin{eqnarray}  \label{kjscfcyhzcfhj}
\ev^{\<1\>}\ =\ \ev\ ,\quad \ev^{\<n+1\>}\ =\ \ev\,(\id\tens \ev^{\<n\>}  \tens \id)\ .
\end{eqnarray}

As in \cite{bbdiff}, we shall use $\square$ for the covariant derivatives on $\mathrm{Vec} A$ and $\Omega^1 A$, reserving $\nabla$ for covariant derivatives on general modules. 
A right covariant derivative
$\square:\Omega^1 A\to \Omega^1 A\tens_A \Omega^1 A$ satisfies the right Leibniz rule, 
\begin{eqnarray}\label{sjgacvaghjchjk}
\square(\xi.a) &=& \square(\xi).a+\xi\tens\extd a\ ,\quad a\in A\ ,\xi\in\Omega^1 A\  .
\end{eqnarray}
A right bimodule covariant derivative $(\square,\sigma^{-1})$ in addition has a  bimodule map
$\sigma^{-1}:\Omega^1 A\tens_A \Omega^1 A\to \Omega^1 A\tens_A \Omega^1 A$ so that 
\begin{eqnarray}\label{sjgacvaghjcdsvhjk}
\square(a.\xi) &=& a.\square(\xi)+\sigma^{-1}(\extd a\tens\xi)\quad \forall\, \xi\in \Omega^1 A,\ a\in A
\end{eqnarray} 
Note that the inverse in $\sigma^{-1}$ is chosen to preserve the conventions of a braided category, even though we do not assume that $\sigma^{-1}$ satisfies the braid relations, nor that it is invertible. 
Now $\square$ extends to a right bimodule covariant derivative $\square^{\<n\>}:\Omega^{\tens n}A\to \Omega^{\tens n+1}A$, defined recursively
\begin{eqnarray} \label{vcakkccjjcj}
\square^{\<0\>} &=& \extd\ ,\quad \square^{\<1\>}\ =\ \square\ ,\cr
\square^{\<n+1\>} &=& \id^{\tens n}\tens\square+(\id^{\tens n}\tens\sigma^{-1})(\square^{\<n\>}\tens\id^{\tens 1})\ .
\end{eqnarray}
The corresponding map $\sigma^{-\<n\>}
:\Omega^{\tens n+1}A\to \Omega^{\tens n+1}A$ is given by the formula
\begin{eqnarray}  \label{bvhkvjhgc}
\sigma^{-\<n\>}\,=\, (\id^{\tens n-1}\tens\sigma^{-1})(\id^{\tens n-2}\tens\sigma^{-1}\tens\id) \dots (\sigma^{-1}\tens \id^{\tens n-1})\ .
\end{eqnarray}
The torsion $\mathrm{Tor}_R$ of $\square$ is defined to be the right $A$-module map
\begin{eqnarray} \label{tordefee}
\mathrm{Tor}_R\ =\ \extd +\wedge\,\square\,:\, \Omega^1A\to \Omega^2 A\ .
\end{eqnarray}

Corresponding to the previous $\square$ on $\Omega^1 A$ (\ref{sjgacvaghjchjk}), there is a left bimodule covariant derivative $\square:\mathrm{Vec} A\to \Omega^1 A\tens_A \mathrm{Vec} A$
so that
\begin{eqnarray} \label{ghjksvdgwv}
\extd\circ\ev\ =\ (\id\tens\ev)(\square\tens\id)+(\ev\tens\id)(\id\tens\square)\ :\ 
\mathrm{Vec} A\tens_A \Omega^1 A\to\Omega^1  A\ .
\end{eqnarray}
This obeys (\ref{vbchiuwjkgc}), where the first line is the left Leibnitz rule, and the second  line, 
where $\sigma:\mathrm{Vec}  A\tens_A \Omega^1 A\to\Omega^1  A\tens_A \mathrm{Vec}  A$ is defined in terms of $\sigma^{-1}$ in (\ref{sjgacvaghjcdsvhjk}) by
\begin{eqnarray}\label{sigrefdefdual}
\sigma\ =\ (\ev\tens\id\tens\id)(\id\tens\sigma^{-1}\tens\id)(\id\tens\id\tens\coev(1))\ ,
\end{eqnarray}
 gives a bimodule covariant derivative:
\begin{eqnarray}\label{vbchiuwjkgc}
\square(a.v) &=& a.\square(v)+\extd a\tens\extd v\ ,\cr
\square(v.a) &=& \square(v).a+\sigma(v\tens\extd a)\quad \forall\, v\in\mathrm{Vec} A,\ a\in A\ .
\end{eqnarray}
This extends to a left bimodule covariant derivative $\square^{\<n\>}:\mathrm{Vec}^{\tens n}A\to \Omega^{1}A\tens_A \mathrm{Vec}^{\tens n}A$, defined recursively by
\begin{eqnarray}  \label{kjvycyjjhgx}
\square^{\<0\>} &=& \extd\ ,\quad \square^{\<1\>}\ =\ \square\ ,\cr
\square^{\<n+1\>} &=& \square\tens \id^{\tens n}+
(\sigma\tens\id^{\tens n})(\id^{\tens 1}\tens\square^{\<n\>})\ .
\end{eqnarray}
The corresponding $\sigma^{\<n\>}:\mathrm{Vec}^{\tens n}A\tens\Omega^1 A\to \Omega^{1}A\tens_A \mathrm{Vec}^{\tens n}A$  is given by
\begin{eqnarray}
(\ev^{\<n\>}\tens\id)(\id^{\tens n}\tens \sigma^{-\<n\>}) &=& (\id\tens\ev^{\<n\>})
(\sigma^{\<n\>}\tens \id^{\tens n})\cr && :\mathrm{Vec}^{\tens n}A\tens_A\Omega^1 A\tens_A\Omega^{\tens n}A\to \Omega^1 A\ .
\end{eqnarray}

Suppose that we wish to take multiple derivatives of a section of a vector bundle. It comes as no surprise that we have to use a covariant derivative $\nabla$ on sections of the vector bundle, but it is quite likely that the result would be written in terms of local coordinates on the manifold. To get away from these local coordinates, we would have to write the derivatives in terms of 1-forms, and then for successive derivatives we would have to use a covariant derivative $\square$ on the 1-forms, as follows.
Suppose that $E$ is a left $A$ module, with a left covariant derivative $\nabla_E:E\to\Omega^1 A\tens_A E$ (i.e.\ $\nabla_E$ satisfies the left Leibnitz rule). 
We iterate this to define $\nabla_E^{(n)}:E\to \Omega^{\tens n}A\tens_A E$
recursively by
\begin{eqnarray}\label{jktycfxtu}
\nabla_E^{(1)}\ =\ \nabla_E\ ,\quad \nabla_E^{(n+1)}=(\square^{\<n\>}\tens\id_E+\id^{\tens n}\tens\nabla_E)\,
\nabla_E^{(n)}\ .
\end{eqnarray}

\begin{defin}    \label{cvgjfxtzhc}
The category ${}_A\mathcal{E}$ consists of objects $(E,\nabla_E)$, where
$E$ is a left  $A$-module, and $\nabla_E$ is a left covariant derivative on $E$. The morphisms $T:(E,\nabla_E)\to (F,\nabla_F)$ consist of left module maps
$T:E\to F$ for which $(\id\tens T)\nabla_E=\nabla_F\,T:E\to \Omega^1 A\tens_A F$.

The category ${}_A\mathcal{E}_A$ consists of objects $(E,\nabla_E,\sigma_E)$, where
$E$ is an  $A$-bimodule, and $(\nabla_E,\sigma_E)$ is a left bimodule covariant derivative on $E$. The morphisms $T:(E,\nabla_E,\sigma_E)\to (F,\nabla_F,\sigma_F)$ consist of bimodule maps
$T:E\to F$ for which $(\id\tens T)\nabla_E=\nabla_F\,T:E\to \Omega^1 A\tens_A F$. It is then automatically true that $\sigma_F(T\tens\id)=(\id\tens T)\sigma_E$. Taking the identity object as 
the algebra $A$ itself, with $\nabla=\extd:A\to \Omega^1 A\tens_A A\cong\Omega^1 A$ and the following tensor product, makes ${}_A\mathcal{E}_A$ into a monoidal category:
\begin{eqnarray*}
\nabla_{E\tens F} &=& \nabla_E\tens
\id_{F}+(\sigma_E\tens\id_{F})(\id_{E}\tens\nabla_F)\ ,\cr
\sigma_{E\tens F} &=& (\sigma_{E}\tens
\id)(\id\tens \sigma_{F})\ .
\end{eqnarray*}
The map $\sigma_E^{-1}:\mathrm{Vec}A\tens_A E\to E\tens_A \mathrm{Vec}A$ is defined by
\begin{eqnarray*}
(\id\tens\ev)(\sigma_E^{-1}\tens\id)\,=\,(\ev\tens\id)(\id\tens\sigma_E):\mathrm{Vec}A\tens_A E\tens_A \Omega^1 A\to E\ .
\end{eqnarray*}

\end{defin}

\section{Noncommutative differential operators}  \label{jkhsvcagkvc}
Here will set out the construction of noncommutative differential operators in \cite{bbdiff}. 
From (\ref{jktycfxtu}) we can define an `action' of $\underline v\in \mathrm{Vec}^{\tens n}A$ on $e\in E$ by
\begin{eqnarray}  \label{bcahkjlclkv}
\underline v\,\la\, e \ =\ (\ev^{\<n\>}\tens\id_E)\, (\underline v\tens\nabla_E^{(n)}e)\ .
\end{eqnarray}
One of the main results of \cite{bbdiff} is that (\ref{bcahkjlclkv}) really is an action of an algebra
\begin{eqnarray*}
\mathcal{T}\,\mathrm{Vec} A\ =\ \bigoplus_{n\ge 0} \mathrm{Vec}^{\tens n}A\ ,
\end{eqnarray*}
but instead of the usual $\tens_A$ product we use an associative product $\bullet$ involving differentiation. This is defined in Lemma \ref{hjsvhjkvjhk} and Theorem \ref{kuykjvcvu}, but the seemingly complicated definition is derived from the principle that (\ref{bcahkjlclkv}) should be an action.

\begin{lemma}\cite{bbdiff} \label{hjsvhjkvjhk}
For all $k\ge 0$, the following recursive procedure
gives a well defined function  $\bullet_k:\mathrm{Vec}^{\tens n}A\tens \mathrm{Vec}^{\tens m}A\to \mathrm{Vec}^{\tens k}A$
satisfying $(a.\underline v)\bullet_k\underline w=a.(\underline v\bullet_k\underline w)$, for all $a\in A$. The definition is recursive in $n\ge 0$: The starting cases are (for $u\in \mathrm{Vec} A$
and $\underline w\in\mathrm{Vec}^{\tens m}A$)
\begin{eqnarray*}
n=0\ , & a\bullet_k \underline w &= \left\{\begin{array}{cc}a.\underline w & k=m \\0 & k\neq m\end{array}\right. \cr
n=1\ , & u\bullet_k \underline w &= \left\{\begin{array}{cc}u\tens w & k=m+1 \\
(\ev\tens\id^{\tens m})(u\tens\square^{\<m\>}\underline w) & k= m
\\0 & \mathrm{otherwise}\end{array}\right. \ .
\end{eqnarray*}
The definition continues with, for $\underline v\in \mathrm{Vec}^{\tens n}A$
(setting $\bullet_{-1}$ to be zero), 
\begin{eqnarray*}
(u\tens\underline v)\bullet_k\underline w &=& u\tens (\underline v\bullet_{k-1}\underline w)
+u\bullet_k(\underline v\bullet_k\underline w)
-(u\bullet_n \underline v)\bullet_k\underline w\ .
\end{eqnarray*}
\end{lemma}

\medskip We modify the usual $A$-bimodule structure on $\mathcal{T}\,\mathrm{Vec} A$ to give
$\mathcal{T}\,\mathrm{Vec} A_\bullet$, where we use the $\bullet$ product for the right $A$-action.

\begin{theorem} \cite{bbdiff} \label{kuykjvcvu}
The $A$-bimodule $\mathcal{T}\,\mathrm{Vec} A_\bullet$ with
product $\bullet:\mathcal{T}\,\mathrm{Vec} A_\bullet\tens_A \mathcal{T}\,\mathrm{Vec} A_\bullet\to \mathcal{T}\,\mathrm{Vec} A_\bullet$
defined by 
\begin{eqnarray*}
\underline v\bullet\underline w \ =\  \sum_{k\ge 0} \underline v\bullet_k\underline w
\end{eqnarray*}
is an associative algebra, with unit $1\in \mathrm{Vec}^{\tens 0}A=A$. Further, for a left $A$-module $E$ with left covariant derivative $\nabla_E$, the map 
in (\ref{bcahkjlclkv}) gives $\la:\mathcal{T}\,\mathrm{Vec} A_\bullet\tens_A E\to E$
which is an action of this algebra. 
\end{theorem}

\medskip In fact $\mathcal{T}\,\mathrm{Vec} A_\bullet$ can be given a left covariant derivative
by the formula
\begin{eqnarray} \label{kjhsacvkyavkuxc}
\nabla(\underline{v})\,=\,\coev(1)\bullet \underline{v}\ ,
\end{eqnarray}
and as this is a right $A$-module map, we see that $(\mathcal{T}\,\mathrm{Vec} A_\bullet,\nabla,0)$
is an object of the category ${}_A\mathcal{E}_A$.  In \cite{bbdiff} a natural transformation $\vartheta_E: \mathcal{T}\mathrm{Vec} A_\bullet\tens_A E\to  E\tens_A\mathcal{T}\mathrm{Vec} A_\bullet$ is defined which makes $(\mathcal{T}\,\mathrm{Vec} A_\bullet,\nabla,0)$ into a central object in  ${}_A\mathcal{E}_A$. Not surprisingly, $\vartheta_E$ is constructed by recursion on $n$ as a map $:\mathrm{Vec}^{\tens n} A\tens E\to
E\tens_A \mathcal{T}\mathrm{Vec} A$.
To do this we start with $n=0$ and $\vartheta_E:A\tens_A E\to E\tens_A A$ being the identity.
For $n=1$,
\begin{eqnarray}  \label{cnjdibvb}
\vartheta_E\ =\ \la+ \sigma_E^{-1}:\mathrm{Vec} A\tens E\to
E\tens_A \mathcal{T}\mathrm{Vec} A\ .
\end{eqnarray}
The definition contunies by, for $u\in \mathrm{Vec} A$,
\begin{eqnarray}  \label{fcauyxaeffrt}
\vartheta_E(w\bullet \underline{v}\tens e) &=& (\la\tens\id)(\id\tens\vartheta_E)(w\tens \underline{v}\tens e) \cr
&&+\ (\sigma_E^{-1}\bullet\id)(\id\tens\vartheta_E)(w\tens \underline{v}\tens e) \ .
\end{eqnarray}
We have the following properties for $\vartheta$:
\begin{eqnarray} \label{jkhzdcvjxty}
\vartheta_E(\underline{u}\bullet \underline{v}\tens e) &=& (\id_E\tens\bullet)(\vartheta_E\tens\id)(\id\tens \vartheta_E)(\underline{u}\tens \underline{v}\tens e)\ ,\cr
\vartheta_{E\tens F} &=& (\id_E\tens \vartheta_{F})(\vartheta_{E}\tens\id_F):
\mathcal{T}\mathrm{Vec} A_\bullet\tens_A E\tens_A F\to E\tens_A F\tens_A \mathcal{T}\mathrm{Vec} A_\bullet\ ,\cr
\la_{E\tens F}  &=& (\id_E\tens\la_F)(\vartheta_E\tens\id_F)  : \mathcal{T}\mathrm{Vec} A_\bullet\tens_A E\tens_A F\to E\tens_A F\ .
\end{eqnarray}

\section{Curvature as a differential operator}
 Given 
the dual basis $\coev(1)=\xi_i\tens u_i\in\Omega^1 A\tens_A\mathrm{Vec}A$ (summation over $i$), the formula for $\nabla:\mathcal{T}\mathrm{Vec}A\to \Omega^1 A\tens_A\mathcal{T}\mathrm{Vec}A$ 
given in (\ref{kjhsacvkyavkuxc}) is
\begin{eqnarray}
\nabla(\underline{v}) &=& \xi_i\tens (u_i\bullet\underline{v})\ .
\end{eqnarray}
The curvature $R:\mathcal{T}\mathrm{Vec}A_\bullet\to\Omega^2 A\tens_A \mathcal{T}\mathrm{Vec}A_\bullet$ is then
(summing over $i,j$)
\begin{eqnarray}
R(\underline{v}) &=& \extd\xi_i\tens (u_i\bullet\underline{v}) -  \xi_i\wedge\nabla (u_i\bullet\underline{v})\cr
&=& \extd\xi_i\tens (u_i\bullet\underline{v}) -  \xi_i\wedge \xi_j\tens (u_j\bullet u_i\bullet\underline{v})\ .
\end{eqnarray}
We can write the curvature, using the associativity of $\bullet$, as
\begin{eqnarray}
R(\underline{v}) &=&\mathcal{R}\bullet  \underline{v}\ ,
\end{eqnarray}
where we have
\begin{eqnarray}\label{vcudisyiuytdstry}
\mathcal{R} &=& \extd\xi_i\tens u_i -  \xi_i\wedge \xi_j\tens u_j\bullet u_i\ \in\ \Omega^2 A\tens_A \mathcal{T}\mathrm{Vec}A_\bullet \ .
\end{eqnarray}
Now, from the formula for $\bullet$,
\begin{eqnarray}
\xi_i\tens \xi_j\tens u_j\bullet u_i &=&\xi_i\tens  \xi_j\tens u_j\tens u_i  + \xi_i\tens \xi_j\tens (\ev\tens\id)(u_j\tens\square u_i ) \cr
&=& \xi_i\tens  \xi_j\tens u_j\tens u_i +(\id^{\tens 2}\tens\ev\tens\id)(\id\tens\coev\tens\id^{\tens 2})(\xi_i\tens \square u_i ) \cr
&=&\xi_i\tens  \xi_j\tens u_j\tens u_i + \xi_i\tens\square u_i  \cr
&=& \xi_i\tens  \xi_j\tens u_j\tens u_i -\square \xi_i\tens u_i \ ,
\end{eqnarray}
where we have used the usual equations for the evaluation and coevaluation. Now we can rewrite $\mathcal{R}$ from (\ref{vcudisyiuytdstry}), using the torsion from (\ref{tordefee}), as
\begin{eqnarray} \label{jkghacfxkytdjyt}
\mathcal{R} &=& \extd\xi_i\tens u_i - \xi_i\wedge  \xi_j\tens (u_j\tens u_i) +\wedge\square\, \xi_i\tens u_i \cr
&=&  \mathrm{Tor}_R(\xi_i)\tens u_i - \xi_i\wedge  \xi_j\tens (u_j\tens u_i) \ .
\end{eqnarray}
Now remember that $R$ is the curvature of a left covariant derivative, and so is a left module map. Applying $R$ to $a\in A$ we get
\begin{eqnarray}\label{juhdtyhgrs}
\mathcal{R}\bullet a \ =\ R(a) \ =\  R(a.1)\ =\ a.R(1)\ =\ a.\mathcal{R}\bullet 1\ =\ a.\mathcal{R}\ ,\end{eqnarray}
so we see that $\mathcal{R}\in \Omega^2 A\tens_A \mathcal{T}\mathrm{Vec}A_\bullet$ is central.

\begin{propos} \label{cvakusvkjghkgv}
Suppose that $(E,\nabla_E)$ is a left $A$-module with left covariant derivative, with curvature $R_E$. Then we can make the $ \mathcal{T}\mathrm{Vec}A_\bullet$ component of $\mathcal{R}$ act on $E$ in the usual manner, giving $\mathcal{R}\,\la\,e=R_E(e)$. 
\end{propos}
\noindent {\bf Proof:}\quad 
We have,  for $e\in E$,
\begin{eqnarray} \label{vuauyvvbil}
\mathcal{R}\,\la\,e &=& \mathrm{Tor}_R(\xi_i)\tens u_i\,\la\,e - \xi_i\wedge  \xi_j\tens (u_j\tens u_i)\,\la\,e\ \in\ \Omega^2 A\tens_A E\ .
\end{eqnarray}
By using the definition of $\la$,
\begin{eqnarray}
\xi_i\wedge  \xi_j\tens (u_j\tens u_i)\,\la\,e &=& \xi_i\wedge  \xi_j\tens (\ev\tens\id_E)(u_j\tens  (\ev\tens\id\tens\id_E)(u_i\tens\nabla_E^{(2)}(e))\ ,\cr
\mathrm{Tor}_R(\xi_i)\tens u_i\,\la\,e &=& \mathrm{Tor}_R(\xi_i)\tens (\ev\tens\id_E)(u_i\tens\nabla_E(e))\ .
\end{eqnarray}
If we set $\nabla_E^{(2)}(e)=\eta_1\tens\eta_2\tens f\in \Omega^{\tens 2}A\tens_A E$
and $\nabla_E(e)=\eta_3\tens g$
(summation implicit), then (using the fact that $\mathrm{Tor}_R$ is a right module map,
\begin{eqnarray}
\xi_i\wedge  \xi_j\tens (u_j\tens u_i)\,\la\,e &=& \xi_i\wedge  \xi_j\tens (\ev\tens\id_E)(u_j\tens  (\ev\tens\id\tens\id_E)(u_i\tens\nabla_E^{(2)}(e))\ ,\cr
&=& \xi_i\wedge  \xi_j\tens u_j(u_i(\eta_1).\eta_2).f \cr
&=& \xi_i\wedge  \xi_j.u_j(u_i(\eta_1).\eta_2)\tens f \cr
&=& \xi_i\wedge  u_i(\eta_1).\eta_2 \tens f \cr
&=& \xi_i.u_i(\eta_1)\wedge  \eta_2 \tens f \cr
&=& \eta_1\wedge  \eta_2 \tens f \ =\ \nabla_E^{(2)}(e)\ ,\cr
\mathrm{Tor}_R(\xi_i)\tens u_i\,\la\,e &=& \mathrm{Tor}_R(\xi_i)\tens u_i(\eta_3).g \cr
 &=& \mathrm{Tor}_R(\xi_i.u_i(\eta_3))\tens g \cr
  &=& \mathrm{Tor}_R(\eta_3)\tens g \ =\ (\mathrm{Tor}_R\tens\id_E)\nabla_E(e)\ .
\end{eqnarray}
Substituting this into (\ref{vuauyvvbil}) gives the result.\quad$\blacksquare$

\section{The sheaf of differential operators $\mathcal{D}A$}  \label{bchksvkuv}
Let $(\Omega^2 A)'=\Hom_A(\Omega^2 A,A)$ (the right $A$-module maps), and define $\widehat{\mathcal{R}}:(\Omega^2 A)'\to  \mathcal{T}\mathrm{Vec}A_\bullet$ by $\widehat{\mathcal{R}}(\alpha)=(\alpha\tens\id)\mathcal{R}$. Then $\widehat{\mathcal{R}}$ is a bimodule map, as $\mathcal{R}$ is central (see (\ref{juhdtyhgrs})). 
Define $\mathcal{W}\subset \mathcal{T}\mathrm{Vec}A_\bullet$ to be the 2-sided ideal (for the $\bullet$ product) generated by the image of $\widehat{\mathcal{R}}$, and $\mathcal{D}A$ to be the quotient 
$\mathcal{D}A=\mathcal{T}\mathrm{Vec}A_\bullet/\mathcal{W}$. We use $[\underline{v}]$ to denote the equivalence class in $\mathcal{D}A$ containing $\underline{v}\in \mathcal{T}\mathrm{Vec}A_\bullet$. 

Recall that $\mathcal{T}\mathrm{Vec}A_\bullet$ has a left action, denoted $\la$, on all objects of
${}_A\mathcal{E}$. We wish to consider modules of the quotient $\mathcal{D}A$, and the obvious way to do this is to define $[\underline{v}]\,\la\,e=\underline{v}\,\la\,e$. However then we have to suppose that $\mathcal{W}$ annihilates every $e\in E$, or equivalently that every element of the image of 
$\widehat{\mathcal{R}}$ annihilates every $e\in E$.

\begin{propos} If $(E,\nabla_E)$ has zero curvature, then the formula $[\underline{v}]\,\la\,e=\underline{v}\,\la\,e$ defines an action of $\mathcal{D}A$. Conversely, under the assumption that $\Omega^2 A$ is finitely generated projective as a right $A$-module, if the formula $[\underline{v}]\,\la\,e=\underline{v}\,\la\,e$ defines an action of $\mathcal{D}A$ on $(E,\nabla_E)$, then $(E,\nabla_E)$ has zero curvature. 
\end{propos}
\noindent {\bf Proof:}\quad 
From Proposition \ref{cvakusvkjghkgv} we have, for $\alpha\in(\Omega^2 A)'$
\begin{eqnarray}
(\alpha\tens\id)\mathcal{R}\,\la\,e\ =\ (\alpha\tens\id_E)R_E(e)\ \in\ E\ ,
\end{eqnarray}
so if $R_E=0$ then the formula gives an action of $\mathcal{D}A$. 

Conversely, if the finitely generated projective assumption holds and $(\alpha\tens\id)R_E(e)=0$ for all $\alpha\in(\Omega^2 A)'$, then $R_E=0$. \quad$\blacksquare$

\medskip 
If we remember that the left covariant derivative $\nabla$ on $\mathcal{T}\mathrm{Vec}A_\bullet$ from (\ref{kjhsacvkyavkuxc}) is
\begin{eqnarray}
\nabla(\underline{v}) &=& \xi_i\tens u_i\bullet \underline{v}\ ,
\end{eqnarray}
it is obvious that $\nabla$ restricts to the ideal $\mathcal{W}\subset \mathcal{T}\mathrm{Vec}A_\bullet$, and has a quotient $\nabla_{\mathcal{D}A}$ on $\mathcal{D}A$. (In fact, we could use any ideal for the $\bullet$ product here.) Explicitly, the covariant derivative is given by
\begin{eqnarray}
\nabla_{\mathcal{D}A}([\underline{v}]) &=& \xi_i\tens [u_i\bullet \underline{v}]\ .
\end{eqnarray}

\begin{propos}
Suppose that $\Omega^2 A$ is finitely generated projective as a right $A$-module. Then
the curvature $R_\mathcal{D}A$ of $(\mathcal{D}A,\nabla_{\mathcal{D}A})$ is zero.
\end{propos}
\noindent {\bf Proof:}\quad The curvature is given by
\begin{eqnarray}
R_{\mathcal{D}A}([\underline{v}]) &=& (\extd\tens\id_\mathcal{D}-\id\wedge\nabla_\mathcal{D})
\nabla_\mathcal{D}([\underline{v}]) \cr
&=&  (\extd\tens\id_\mathcal{D}-\id\wedge\nabla_\mathcal{D})
\big( \xi_i\tens [u_i\bullet \underline{v}]\big) \cr
&=& \extd \xi_i\tens [u_i\bullet \underline{v}] - \xi_i\wedge\xi_j\tens [u_j\bullet u_i\bullet \underline{v}] \cr
&=& \big(\extd \xi_i\tens [u_i] - \xi_i\wedge\xi_j\tens [u_j\bullet u_i] \big) \bullet [\underline{v}] \ .
\end{eqnarray}
From this, it is enough to show that $\mathcal{R}\in \Omega^2 A\tens_A \mathcal{W}$. If we write the dual basis for $\Omega^2 A$ as $\phi_k\tens \alpha_k\in \Omega^2 A\tens_A(\Omega^2 A)'$ (summed over $k$), then
\begin{eqnarray*}
\mathcal{R} &=& \phi_k\tens (\alpha_k\tens\id)\mathcal{R}\ =\ \phi_k\tens 
\widehat{\mathcal{R}}(\alpha_k)\ \in\ \Omega^2 A\tens_A \mathcal{W}\ .\quad\blacksquare
\end{eqnarray*}

\begin{example}{\textbf{Classical differential operators.}} As $\mathcal{R}$ is defined in terms of the curvature of a connection, it is independent of the choice of dual basis. On an open subset of a differential manifold we can choose coordinates $x_1,\dots,x_n$, and dual basis $\xi_i\tens u_i=\extd x_i\tens \frac{\partial}{\partial x_i}$, and then from (\ref{vcudisyiuytdstry}) we have
\begin{eqnarray}
\mathcal{R}\ =\ -\ \extd x_i\wedge\extd x_j\tens \frac{\partial}{\partial x_j} \bullet \frac{\partial}{\partial x_i}\ .
\end{eqnarray}
Now using the basis $\extd x_i\wedge\extd x_j$ for $i<j$ of the two forms and using this to calculate 
$\widehat{\mathcal{R}}$ shows that $\mathcal{W}$ is generated by the following relations,  not entirely unexpectedly: 
\begin{eqnarray}
\frac{\partial}{\partial x_j} \bullet \frac{\partial}{\partial x_i} - \frac{\partial}{\partial x_i} \bullet \frac{\partial}{\partial x_j}\ =\ 0\ .
\end{eqnarray}
Of course, this becomes rather more complicated once functions of the $x_i$ are used to multiply the 
$\frac{\partial}{\partial x_j}$, and the resulting algebra $\mathcal{D}A$ is noncommutative. 
\end{example}

We want to get a rather easier picture of the relations in the noncommutative case, and to do that we refer to the generalisation of Lie bracket of vector fields given in \cite{begexp}.

\begin{defin} \label{antidef}\cite{begexp} An $x\in\mathrm{Vec}A\tens\mathrm{Vec}A$ is called
antisymmetric if $\ev^{\<2\>}(\pi x\tens k)=0$ for all
$k\in \ker\wedge:\Omega^{1}A\tens_A\Omega^{1}A\to\Omega^{2}A$, where $\pi$ is the quotient map from $\mathrm{Vec}A\tens\mathrm{Vec}A$ to $\mathrm{Vec}A\tens_{A}\mathrm{Vec}A$.  We call $\mathrm{AntiVec}_2A$ the set of
antisymmetric elements in $\mathrm{Vec}A\tens\mathrm{Vec}A$.
\end{defin}

\begin{defin}\label{deflie}\cite{begexp} Define  $\varphi:\mathrm{AntiVec}_2A\to {\rm Vec}\,A$
by the following formula,
\[
\varphi(u\tens v)(\xi)\,=\,D_u(v(\xi))\,+\,\ev(\id\tens \ev\tens\id) (u\tens
v\tens z)\ ,\quad\xi\in\Omega^{1}A\ .
\]
where $z\in \Omega^{1}A
\tens_A\Omega^{1}A$ is chosen so that $\wedge z=d\xi$ (the choice does not matter). We use the directional derivative
$D_u(a)=u(\extd a)$.
 To check that its image is in $\mathrm{Vec}A$ we use the
following proposition.
\end{defin}

\begin{propos}\label{prop17} \cite{begexp} The image of the map $\varphi$ in Definition \ref{deflie} is in
    $\mathrm{Vec}A$. Further $\varphi$ is a left
$A$-module map, but not in general a right module map,
as $\varphi(u\tens v).a=\varphi(u\tens v.a)+u.D_v(a)$. 
Also $\varphi(u\tens a.v)=\varphi(u.a\tens v)+D_{u}(a).v$.
\end{propos}
\proof To see that $\varphi(u\tens v)$ is a right module map use the following,
where $\wedge z=\extd\xi$,
\begin{eqnarray*} 
    \varphi(u\tens v)(\xi.a) &=& D_u(v(\xi).a)\,+\,\ev(\id\tens
\ev\tens\id) (u\tens v\tens (z.a-\xi\tens \extd a)) \cr
&=& \varphi(u\tens v)(\xi).a\,+\, u(v(\xi).\extd a)\,-\, u(v(\xi).\extd a)\ .
\end{eqnarray*}
It is quite easy to see that $\varphi(a.u\tens v)(\xi)=a.\varphi(u\tens v)(\xi)$. 
For the right action,
\begin{eqnarray*} 
\varphi(u\tens v)(a.\xi) &=& D_u(v(a.\xi))\,+\,\ev(\id\tens \ev\tens\id) (u\tens
v\tens (a.z+\extd a\tens\xi)) \cr
&=& \varphi(u\tens v.a)(\xi)\,+\,\ev(\id\tens \ev\tens\id) (u\tens
v\tens \extd a\tens\xi)\ .
\end{eqnarray*}
Finally we calculate
\begin{eqnarray*} 
\varphi(u\tens a.v)(\xi) &=& 
D_u(a.v(\xi))\,+\,\ev(\id\tens \ev\tens\id) (u\tens
a.v\tens z) \cr &=& D_u(a).v(\xi)+
D_{u.a}(v(\xi))\,+\,\ev(\id\tens \ev\tens\id) (u.a\tens v\tens z)\
.\quad\blacksquare
\end{eqnarray*}

\medskip Classically, $[u,v]=\varphi(u\tens v-v\tens u)$. We can choose the 1-form $\xi$ to be locally constant in some coordinate system, and then from Definition \ref{deflie} we get the usual result $[u,v](\xi)=D_u(v(\xi))-D_v(u(\xi))$. Note that the major difference in the noncommutative case is that we have not, in general, got an antisymmetrisation procedure. However we can characterise antisymmetric tensor products as follows:

\begin{propos}  \label{jacxcjtyhxz}
Suppose that $\Omega^1 A$ and $\Omega^2 A$ are finitely generated projective as right $A$-modules, with a dual basis $\xi_i\tens u_i\in \Omega^1 A\tens_A \mathrm{Vec}A$. 
Then the subset of all $x$ in $\mathrm{Vec}A\tens_{A}\mathrm{Vec}A$ for which
$\ev^{\<2\>}(x\tens k)=0$ for all $k\in \ker\wedge:\Omega^1 A\tens_A \Omega^1 A\to \Omega^2 A$ is
\begin{eqnarray*}
\big\{\alpha(\xi_i\wedge  \xi_j)\, u_j\tens u_i : \alpha\in\Hom_A(\Omega^2 A,A)\big\}\ .
\end{eqnarray*}
\end{propos}
\noindent {\bf Proof:}\quad We have a short exact sequence of right modules, 
\begin{eqnarray*}
0 \longrightarrow \ker\wedge \longrightarrow \Omega^1 A\tens_A \Omega^1 A \stackrel{\wedge}
\longrightarrow \Omega^2 A\longrightarrow 0\ ,
\end{eqnarray*}
where $\Omega^2 A$ is finitely generated projective as a right $A$-module. We deduce that there is a splitting map $S:\Omega^2 A\to \Omega^1 A\tens_A \Omega^1 A$. Given $x$ so that $\ev^{\<2\>}(x\tens k)=0$ for all $k\in \ker\wedge$, we get a well defined map $\omega\mapsto \ev^{\<2\>}(x\tens S(\omega))$ in $\Hom_A(\Omega^2 A,A)$. The map $\alpha\in \Hom_A(\Omega^2 A,A)$ arises in this way from $x=\alpha(\xi_i\wedge  \xi_j)\, u_j\tens u_i$. Any $x$ corresponding to $\alpha=0$ would have to be zero on all of $\Omega^1 A\tens_A \Omega^1 A$, i.e.\ $x=0$.\quad$\blacksquare$

\begin{propos} \label{jkghabvdcyct}     For $\alpha\in\Hom_A(\Omega^2 A,A)$,
\begin{eqnarray*}
(\alpha\tens\id)\mathcal{R} &=& \varphi(\alpha(\xi_i\wedge  \xi_j)\, u_j\tens u_i) -  \alpha(\xi_i\wedge  \xi_j)\, u_j\bullet u_i 
\end{eqnarray*}
\end{propos}
\noindent {\bf Proof:}\quad 
Given $\alpha\in\Hom_A(\Omega^2 A,A)$, from (\ref{jkghacfxkytdjyt}) we have 
\begin{eqnarray} \label{jkghacyct}
(\alpha\tens\id)\mathcal{R} &=& \alpha(\mathrm{Tor}_R(\xi_i))\, u_i - \alpha(\xi_i\wedge  \xi_j)\, u_j\tens u_i  \ .
\end{eqnarray}
By the usual evaluation and coevaluation properties, $\alpha(\xi_i\wedge  \xi_j)\, u_j\tens u_i\in 
\mathrm{AntiVec}_2A$ (see Definition \ref{antidef}). Then, for $\eta\in\Omega^1 A$,
and $\wedge z=\extd\eta$,
\begin{eqnarray}
\varphi(\alpha(\xi_i\wedge  \xi_j)\, u_j\tens u_i)(\eta) &=& \alpha(\xi_i\wedge  \xi_j)\,\big(
D_{u_j}(u_i(\eta))+\ev^{\<2\>}(u_j\tens u_i\tens z)\big) \cr
&=& \alpha(\xi_i\wedge  \xi_j)\,D_{u_j}(u_i(\eta))+\alpha(\extd\eta)\ .
\end{eqnarray}
We expand the $D_{u_j}(u_i(\eta))$ term using $\square$,
\begin{eqnarray}
\varphi(\alpha(\xi_i\wedge  \xi_j)\, u_j\tens u_i)(\eta) 
&=& \alpha(\xi_i\wedge  \xi_j)\,\big((\ev\tens\ev)(u_j\tens \square u_i\tens\eta)
\cr &&+\ev^{\<2\>}(u_j\tens  u_i\tens \square\eta)
\big)+\alpha(\extd\eta) \cr
&=& \alpha(\xi_i\wedge  \xi_j)\,(\ev\tens\ev)(u_j\tens \square u_i\tens\eta)
\cr &&+\alpha(\wedge \square\eta)+\alpha(\extd\eta) \cr
&=&  \alpha(\xi_i\wedge  \xi_j)\,(\ev\tens\ev)(u_j\tens \square u_i\tens\eta)+\alpha(\mathrm{Tor}_R(\eta))\ .
\end{eqnarray}
This gives
\begin{eqnarray}
\varphi(\alpha(\xi_i\wedge  \xi_j)\, u_j\tens u_i)
&=&  \alpha(\xi_i\wedge  \xi_j)\,(\ev\tens\id)(u_j\tens \square u_i)+\alpha(\mathrm{Tor}_R(\xi_i))\,u_i\ ,
\end{eqnarray}
and on substitution in (\ref{jkghacyct}) we get
\begin{eqnarray*}
(\alpha\tens\id)\mathcal{R} &=& \varphi(\alpha(\xi_i\wedge  \xi_j)\, u_j\tens u_i) - \alpha(\xi_i\wedge  \xi_j)\,(\ev\tens\id)(u_j\tens \square u_i) \cr && -\ \alpha(\xi_i\wedge  \xi_j)\, u_j\tens u_i  \cr
&=& \varphi(\alpha(\xi_i\wedge  \xi_j)\, u_j\tens u_i) -  \alpha(\xi_i\wedge  \xi_j)\, u_j\bullet u_i \ .\quad\blacksquare
\end{eqnarray*}

\begin{cor}
Suppose that $\Omega^1 A$ and $\Omega^2 A$ are finitely generated projective as right $A$-modules. Then
the relations in $\mathcal{D}A$ are of the form $\varphi(u\tens v)-u\bullet v$ for every
$u\tens v\in \mathrm{Vec}A\tens_A \mathrm{Vec}A$ (summation implicit) so that $\ev^{\<2\>}(u\tens v\tens k)=0$ for all $k\in \ker\wedge:\Omega^1 A\tens_A \Omega^1 A\to \Omega^2 A$.
\end{cor}
\noindent {\bf Proof:}\quad Propositions \ref{jkghabvdcyct} and \ref{jacxcjtyhxz}. Note that the relation depends only on $u\tens v\in \mathrm{Vec}A\tens_A \mathrm{Vec}A$, with emphasis on the $\tens_A$. This is as, by Proposition \ref{prop17}, the contributions from the two terms  cancel when applied to $u.a\tens v-u\tens a.v$.  \quad$\blacksquare$

\section{Example: Differential calculi on Hopf algebras}  \label{bhsoabvi}
The idea of these calculi follows from \cite{W2}.
For a Hopf algebra $H$ with a left covariant differential calculus, there is a well defined left coaction $\lambda:\Omega^1 H\to H\tens \Omega^1 H$ defined by $\lambda(x.\extd y)=x_{(1)}\, y_{(1)}\tens x_{(2)}.\extd y_{(2)}$. We write
$\lambda(\xi)=\xi_{[-1]}\tens\xi_{[0]}$. Set $L^1H$ to be the coinvariants under the left coaction, i.e.\ those $\xi\in \Omega^1 H$ so that $\lambda(\xi)=1\tens\xi$. If the antipode $S$ is invertible, there is an isomorphism $\Omega^1 H\cong L^1 H\tens H$, given by product one way, and the other way by
$\xi\mapsto \xi_{[0]}.S^{-1}(\xi_{[-1]})\tens \xi_{[-2]}$. If $\Omega^1 H$ is finitely generated as a right $H$-module, then $L^1 H$ is a finite dimensional vector space. There is a left action of $H$ on $L^1 H$ given by $h\,\la\,\eta=h_{(2)} \,\eta\,S^{-1}(h_{(1)})$. The product of a element of $H$ on the left with an element of $\Omega^1 H\cong L^1 H\tens H$ is given by $h.(\xi\tens g)=h_{(2)}\,\la\,\xi\tens h_{(1)}\,g$. 

If we set $\mathfrak{h}$ to be the vector space dual of $L^1 H$, then $\mathrm{Vec}H\cong H\tens \mathfrak{h}$, where the pairing between $\mathrm{Vec}H$ and $\Omega^1 H$ is given by
$(h\tens\alpha)(\xi\tens g)=h\,\alpha(\xi)\,g\in H$. We take $\xi_i\tens u_i\in L^1 H\tens \mathfrak{h}$ (summing over $i$) to be a dual vector space basis. Now from (\ref{vcudisyiuytdstry})
\begin{eqnarray}\label{vcudisyiuytdstrdsy}
\mathcal{R} &=& \extd\xi_i\tens u_i -  \xi_i\wedge \xi_j\tens u_j\bullet u_i\ .
\end{eqnarray}
Now take a vector space basis $\omega^k$ of the left invariant 2-forms, and
(summing over $k$) 
\begin{eqnarray}  \label{bciuasdog}
\extd\xi_i\ =\ n_{ik}\,\omega_k\ ,\quad \xi_i\wedge\xi_j\ =\ e_{ijk}\,\omega_k\ .
\end{eqnarray}
The relations in $\mathcal{D}A$ are then
\begin{eqnarray}  
n_{ik}\,u_i\ =\ e_{ijk}\,u_j\bullet u_i\ .
\end{eqnarray}
Further analysis of this for functions on a Lie group would give the classical Lie algebra (i.e.\ $\mathfrak{h}$) 
structure constants. 

For a noncommutative example, we will give the relations for the 3D calculus on quantum $SU_2$ (see \cite{W2} for the original work, and \cite{BMriem} for our notation). There is a basis $e^0,e^+,e^-$ for the left invariant 1-forms, and we then have:
\[ \extd  e^0=q^3 e^+\wedge e^-,\quad
\extd e^\pm=\mp\, q^{\pm 2}\,(1+q^{- 2})\,e^\pm\wedge e^0,\quad
e^\pm\wedge e^\pm= e^0\wedge e^0=0\]
\[
q^2 e^+\wedge e^-+ e^-\wedge e^+=0,\quad  e^0\wedge e ^\pm+q^{\pm4}\,e^\pm\wedge e^0=0
\] 
We now set $\omega_0=e^+\wedge e^-$ and $\omega_\pm=e^\pm\wedge e^0$, and take $u_0,u_\pm$ to be the dual basis of $e^0,e^+,e^-$ (i.e.\ $u_i(e^j)=\delta_{ij}$). Now the coefficients in (\ref{bciuasdog}) are given by $e_{iik}=0$ and
\begin{eqnarray*}
n_{0k}\ =\ \left\{\begin{array}{cc}q^3 & k=0 \\0 & \mathrm{otherwise}\end{array}\right.\ ,\quad
n_{\pm k}\ =\ \left\{\begin{array}{cc} \mp\, q^{\pm 2}\,(1+q^{- 2}) & k=\pm \\0 & \mathrm{otherwise}\end{array}\right.\ ,
\end{eqnarray*}
\begin{eqnarray*}
e_{\pm 0 k}\ =\ \left\{\begin{array}{cc}1 & k=\pm \\0 & \mathrm{otherwise}\end{array}\right.\ ,\quad
e_{0\pm  k}\ =\ \left\{\begin{array}{cc}-q^{\pm4} & k=\pm \\0 & \mathrm{otherwise}\end{array}\right.\ ,
\end{eqnarray*}
\begin{eqnarray*}
e_{+-k}\ =\ \left\{\begin{array}{cc}1 & k=0 \\0 & \mathrm{otherwise}\end{array}\right.\ ,\quad
e_{-+k}\ =\ \left\{\begin{array}{cc}-q^2 & k=0 \\0 & \mathrm{otherwise}\end{array}\right.\ .
\end{eqnarray*}
This gives the relations in $\mathcal{D}A$,
\begin{eqnarray} \label{kuycvghjkjhc}
q^3\,u_0 &=& u_-\bullet u_+-q^2\,u_+\bullet u_-\ ,\cr
 \mp\, q^{\pm 2}\,(1+q^{- 2})\,u_\pm &=& u_0\bullet u_\pm - q^{\pm4}\,u_\pm\bullet u_0\ .
\end{eqnarray}

In \cite{BMriem} it is shown that the only left invariant left covariant derivatives on the 3D $\Omega^1 \C_q[SU_2]$ which are bimodule covariant derivatives, and are invariant under the right $\mathbb{CZ}$ coaction are of the form (\ref{bchkdlsvkvk}) (where we use $\nabla^L$ for the restriction to the left invariant forms). Note we have set $\nu=\mu_+$, $\mu=\mu_-$ compared to the notation in \cite{BMriem}, and assume that $q$ is not a root of unity. 
\begin{eqnarray}  \label{bchkdlsvkvk}
\nabla^L(e^0)  &=& r\, e^0\tens e^0+ \mu_+\, e^+\tens e^- + \mu_-\,e^-\tens e^+\ ,\cr
\nabla^L(e^\pm) &=& n_\pm\, e^0\tens e^\pm + m_\pm\, e^\pm\tens e^0\ .
\end{eqnarray}

\begin{propos}  \label{bchkdlsvkvk1}
The curvature of the connection in (\ref{bchkdlsvkvk}) is (summing over the $\pm$ in the second formula)
\begin{eqnarray*}
R(e^\pm) 
&=&  \big(  n_\pm\, q^3 e^+\wedge e^-  - \mu_\mp\,m_\pm\, e^\pm\wedge  e^\mp \big)\tens e^\pm
\cr && +\,  m_\pm\, \big( \mp  q^{\pm 2}\,(1+q^{- 2}) + n_\pm \, q^{\pm4} - r  \big) e^\pm\wedge e^0\tens e^0\ ,\cr
R(e^0)  
  &=& \big(r\, q^3  - \mu_+\, m_-
  + \mu_-\, m_+\, q^2 
 \big)e^+\wedge e^-\tens e^0   \cr
 &&+\, \mu_\pm\,\big( \mp \, q^{\pm 2}\,(1+q^{- 2})
 + r\, q^{\pm4} -  n_\mp
 \big)e^\pm \wedge e^0\tens e^\mp\ .
\end{eqnarray*}
\end{propos}
\noindent {\bf Proof:}\quad From the definition of curvature,
\begin{eqnarray*}
R(e^\pm) &=& n_\pm\, \extd e^0\tens e^\pm + m_\pm\, \extd e^\pm\tens e^0\cr
&& -\,n_\pm\, e^0\wedge \nabla^L(e^\pm) - m_\pm\, e^\pm\wedge \nabla^L(e^0)\cr
&=& n_\pm\, q^3 e^+\wedge e^- \tens e^\pm \mp m_\pm\, q^{\pm 2}\,(1+q^{- 2})\,e^\pm\wedge e^0\tens e^0\cr
&& -\,n_\pm\, m_\pm \, e^0\wedge  e^\pm\tens e^0 - m_\pm\, e^\pm\wedge (r\, e^0\tens e^0 + \mu_\mp\,e^\mp\tens e^\pm)\cr
&=&  \big(  n_\pm\, q^3 e^+\wedge e^-  - \mu_\mp\,m_\pm\, e^\pm\wedge  e^\mp \big)\tens e^\pm
\cr && +\,  m_\pm\, \big( \mp  q^{\pm 2}\,(1+q^{- 2})\,e^\pm\wedge e^0 - n_\pm \, e^0\wedge  e^\pm - r\,  e^\pm\wedge e^0\big)\tens e^0\cr
&=&  \big(  n_\pm\, q^3 e^+\wedge e^-  - \mu_\mp\,m_\pm\, e^\pm\wedge  e^\mp \big)\tens e^\pm
\cr && +\,  m_\pm\, \big( \mp  q^{\pm 2}\,(1+q^{- 2}) + n_\pm \, q^{\pm4} - r  \big) e^\pm\wedge e^0\tens e^0\ .
\end{eqnarray*}
Summing over the $\pm$ in the following formula,
\begin{eqnarray*}
R(e^0)  &=& r\, \extd e^0\tens e^0+ \mu_\pm\, \extd e^\pm\tens e^\mp \cr
 && - r\, e^0\wedge \nabla^L( e^0)- \mu_\pm\, e^\pm \wedge \nabla^L( e^\mp) \cr
 &=& r\, q^3 e^+\wedge e^-\tens e^0  \mp \mu_\pm\,\, q^{\pm 2}\,(1+q^{- 2})\,e^\pm\wedge e^0\tens e^\mp \cr
 && - r\, \mu_\pm\, e^0\wedge  e^\pm\tens e^\mp- \mu_\pm\, e^\pm \wedge (n_\mp\, e^0\tens e^\mp + m_\mp\, e^\mp\tens e^0) \cr
 &=& \big(r\, q^3 e^+\wedge e^- - \mu_\pm\, e^\pm \wedge m_\mp\, e^\mp
 \big)\tens e^0   \cr
 &&+\, \mu_\pm\,\big( \mp \, q^{\pm 2}\,(1+q^{- 2})\,e^\pm\wedge e^0
 - r\, e^0\wedge  e^\pm
 -  n_\mp\, e^\pm \wedge e^0
 \big)\tens e^\mp\ .\quad\blacksquare
\end{eqnarray*}

\begin{cor} \label{vcgdsgvjgcx}
The  left covariant derivative in (\ref{bchkdlsvkvk}) has zero curvature in precisely four cases:
\newline a) $\nabla^L=0$. 
\newline b) $n_-\,=\,  - \, q^{2}\,(1+q^{- 2})=-\, \mu_+\,m_-\, q^{-1}$\ ,\ $n_+\,=\,  q^{- 2}\,(1+q^{- 2})=q^{-3} \,\mu_-\,m_+$\ ,\ $r=0$.
\newline c) $n_-=m_-=\mu_+=0$\ ,\ $n_+=q^{-2}=q^{-3} \,\mu_-\,m_+$\ ,\ $r=-1$.
\newline d) $n_+=m_+=\mu_+=0$\ ,\ $n_-= -\,1=-\, \mu_+\,m_-\, q^{-1}$\ ,\ $r=q^{-2}$.
\newline Generically, all these cases have invertible $\sigma$. 
\end{cor}

 From Corollary \ref{vcgdsgvjgcx} and \cite{BMriem}, the only curvature zero case where $\sigma$ satisfies the braid relations is the trivial case $\nabla^L=0$. However we do find the following statement about the star operation $\star:\Omega^1 H\to \overline{\Omega^1 H}$ given by the star operation on quantum $SU_2$, which assumes that $q$ is real.

\begin{cor}   
The condition that the zero curvature left covariant derivative in Corollary \ref{vcgdsgvjgcx} preserves the star operation in the sense that
\begin{eqnarray*}
(\id\tens\star)\,\nabla_{\Omega^1 H}\ =\ \nabla_{\overline{\Omega^1 H}}\ \star
\end{eqnarray*}
is that we have case \textit{(a)} of Corollary \ref{vcgdsgvjgcx}, or the sub-case of \textit{(b)} where 
\begin{eqnarray*}
m_+\,=\,-\, m_-^*\ ,\quad \mu_+\,=\,-\,q^2\,\mu_-^*\ ,\quad -\,\mu_-\,m_-^*\,=\, q+q^{-1}\ .
\end{eqnarray*}
\end{cor}
\noindent {\bf Proof:}\quad Substitute the cases of Corollary \ref{vcgdsgvjgcx} into the equations of \cite{BMriem}.\quad $\blacksquare$

\medskip Now we take the basis order $\{e^+,e^0,e^-\}$ for $L^1 H$, and write matrices to give the action of the $u_i$, by right multiplication on column vectors:
\begin{eqnarray}
u_+ \to \left(\begin{array}{ccc}0 & 0 & 0 \\m_+ & 0 & 0 \\0 & \mu_+ & 0\end{array}\right)\ ,\
u_0 \to \left(\begin{array}{ccc}n_+ & 0 & 0 \\0 & r & 0 \\0 & 0 & n_-\end{array}\right)\ ,\ 
u_-\to \left(\begin{array}{ccc}0 & \mu_- & 0 \\0 & 0 & m_- \\0 & 0 & 0\end{array}\right)\ .
\end{eqnarray}
This is just read off (\ref{bchkdlsvkvk}). However we can now check our calculations by substituting these matrices into the relations in (\ref{kuycvghjkjhc}). This gives the following series of equations, which is precisely the same as the conditions for zero curvature in Corollary \ref{vcgdsgvjgcx}:
\begin{eqnarray}
&& m_+\, \mu_- = n_+ q^3  \ ,\  m_-\, \mu_+ - m_+\, \mu_- q^2 - q^3 r=0  \ , \  m_- \mu_+ =- n_- q   \cr
&&
m_+ (-1 - q^2 + n_+ q^4 - r)=0\ ,\ \mu_+ (1 + n_- + q^2 - q^4 r)=0\ ,\cr
&&
\mu_- (-1 - q^2 + n_+ q^4 - r)=0\ ,\ m_- (1 + n_- + q^2 - q^4 r)=0\ .
\end{eqnarray}

\medskip Why did we only have to look at the $n$-forms to get the relations (\ref{kuycvghjkjhc}), and then have to do lots of hard work to find the curvature? 
The construction of the algebra $\mathcal{D}A$ and its actions was very explicit, and involves specifying covariant derivatives. However it turns out that many of the relations on $\mathcal{D}A$, viewed as an algebra generated by elements of $\mathrm{Vec}A$, do not depend on specifying covariant derivatives. But there are, in general, other relations in $\mathcal{D}A$ -- a choice of dual basis of $\mathrm{Vec}A$ still leaves relations between the basis elements (for a non-freely generated case, the Hopf algebra case is freely generated). We do not claim to have studied the properties of $\mathcal{D}A$ as an abstract algebra, and it might be interesting to do so.

\section{The crossing map for $\mathcal{D}A$} \label{bvhdskv}

\begin{defin}    \label{cvgjfxtzhcflat}
The categories ${}_A\mathcal{F}$ and ${}_A\mathcal{F}_A$ are defined exactly as
${}_A\mathcal{E}$ and ${}_A\mathcal{E}_A$ in Definition \ref{cvgjfxtzhc}, except that we require that the curvature of the left derivative vanishes for both, and for $(E,\nabla_E,\sigma_E)$ in ${}_A\mathcal{F}_A$ we require the existence of a map $\sigma_E:E\tens_A\Omega^2 A\to \Omega^2 A\tens_A E$ so that, on 
$E\tens_A\Omega^1 A\tens_A \Omega^1 A$
\begin{eqnarray} \label{ajhscfxhgfz}
\sigma_E(\id_E\tens \wedge)\ =\ (\wedge\tens\id_E)(\id\tens\sigma_E)(\sigma_E\tens\id)\ .
\end{eqnarray}
\end{defin}

Note that if $\Omega^2 A$ is finitely generated projective as a right $A$-module, then 
(\ref{ajhscfxhgfz}) implies the existence of $\sigma_E^{-1}:(\Omega^2 A)'\tens_A E\to E\tens_A (\Omega^2 A)'$ so that
\begin{eqnarray}  \label{ajhscfxhgfz1}
(\id_E\tens \ev)(\sigma_E^{-1}\tens\id)\,=\, (\ev\tens\id_E)(\id\tens \sigma_E):(\Omega^2 A)'\tens_A E\tens_A \Omega^2 A\to E\ .
\end{eqnarray}

In Section \ref{jkhsvcagkvc}, summarising the results of  \cite{bbdiff}, we stated that 
$(\mathcal{T}\,\mathrm{Vec} A_\bullet,\nabla,0)$, with the left bimodule covariant derivative given in (\ref{kjhsacvkyavkuxc}), was in the center of the category ${}_A\mathcal{E}_A$, using the crossing natural transformation $\vartheta$. In Section \ref{bchksvkuv} we showed (under the condition that $\Omega^2 A$ is finitely generated projective as a right $A$-module) that 
$(\mathcal{D}A,\nabla_\mathcal{D},0)$ is in ${}_A\mathcal{F}_A$. We also showed that $\mathcal{D}A$ acts on objects in ${}_A\mathcal{F}$. 
 It is natural to ask if $(\mathcal{D}A,\nabla_\mathcal{D},0)$ is in the center of ${}_A\mathcal{F}_A$, and to see if the original crossing natural transformation $\vartheta$, defined in (\ref{cnjdibvb}) and (\ref{fcauyxaeffrt}), is compatible with the relations $\widehat{\mathcal{R}}(\alpha)$.

\begin{propos}   \label{bvchas}
Suppose that for $(E,\nabla_E,\sigma_E)\in {}_A\mathcal{E}_A$,
 (\ref{ajhscfxhgfz}) holds, that $\Omega^2 A$ is finitely generated projective as a right $A$-module, and that the curvature $R_E:E\to \Omega^2 A\tens_A E$ is a right module map. Then
\begin{eqnarray*}
\vartheta_E\,(\widehat{\mathcal{R}}\tens \id_E) &=& \widehat{\mathcal{R}}\, \la\, \id_E\tens 1_A
+(\ev\tens\id_E\tens \id)
(\id\tens\sigma_E\tens\widehat{\mathcal{R}})(\id\tens\id_E\tens\coev_2) \cr
&& :\, (\Omega^2 A)'\tens_A E\to E\tens_A \mathcal{T}\mathrm{Vec}A\ .
\end{eqnarray*}
\end{propos}
where $\coev_2(1)\in \Omega^2 A\tens_A (\Omega^2 A)'$ is the dual basis for 
$\Omega^2 A$. 

\noindent {\bf Proof:}\quad 
First, from (\ref{cnjdibvb}) and (\ref{fcauyxaeffrt}), for $\alpha\in (\Omega^2 A)'$ and $e\in E$,
\begin{eqnarray*}
\vartheta_E((\alpha\tens\id)\mathcal{R}\tens e) &=& \alpha(\extd\xi_i)\,\vartheta_E(u_i\tens e)-
\alpha(\xi_i\wedge\xi_j)\,\vartheta_E(u_j\bullet u_i\tens e)   \cr
&=& \alpha(\extd\xi_i)\,\vartheta_E(u_i\tens e)- \cr
&&\alpha(\xi_i\wedge\xi_j)\,(\la\tens\id+\sigma_E^{-1}\bullet\id)(u_j\tens \vartheta_E(u_i\tens e))   \cr
&=& \alpha(\extd\xi_i)\,(u_i\,\la\,e+\sigma_E^{-1}(u_i\tens e))   - \cr
&&\alpha(\xi_i\wedge\xi_j)\,(\la\tens\id+\sigma_E^{-1}\bullet\id)(u_j\tens (u_i\,\la\,e+\sigma_E^{-1}(u_i\tens e)))   \ ,
\end{eqnarray*}
and from this we have
\begin{eqnarray} \label{jhcfshjtfdj}
&& \vartheta_E((\alpha\tens\id)\mathcal{R}\tens e) - (\alpha\tens\id)\mathcal{R}\,\la\,e\tens 1_A\cr
&=& \alpha(\extd\xi_i)\,\sigma_E^{-1}(u_i\tens e)   - \alpha(\xi_i\wedge\xi_j)\,\sigma_E^{-1}(u_j\tens (u_i\,\la\,e)) \,-  \cr
&&\alpha(\xi_i\wedge\xi_j)\,(\la\tens\id+\sigma_E^{-1}\bullet\id)(u_j\tens \sigma_E^{-1}(u_i\tens e))   \ .
\end{eqnarray}
In (\ref{jhcfshjtfdj}) the actions are made of $u\,\la\,e=(\ev\tens\id_E)(u\tens\nabla_Ee)$, and these can be combined with the coevaluations $\coev(1)=\xi_i\tens u_i$ to simplify certain of the terms, as follows. Counting multiplying out the brackets, we have four terms on the RHS of (\ref{jhcfshjtfdj}), and then the second and third terms give
\begin{eqnarray*}
\alpha(\xi_i\wedge\xi_j)\,\sigma_E^{-1}(u_j\tens (u_i\,\la\,e)) &=& 
(\alpha\tens \sigma_E^{-1})(\id\wedge\coev\tens\id_E)\nabla_E(e) \cr
&=& (\alpha\tens\id_E\tens\id)(\id\wedge\sigma_E\tens\id)(\nabla_E(e)\tens\xi_i\tens u_i)
\ ,\cr
\alpha(\xi_i\wedge\xi_j)\,(u_j\,\la\, \sigma_E^{-1}(u_i\tens e)) &=&
(\alpha\tens\id_E\tens\id)(\id\wedge\nabla_E\tens\id) (\xi_i\tens\sigma_E^{-1}(u_i\tens e))\ .
\end{eqnarray*}
The fourth term of (\ref{jhcfshjtfdj}) gives (using the definition of $\bullet$)
\begin{eqnarray}  \label{hvcfhjjhfxxgh}
&& \alpha(\xi_i\wedge\xi_j)\,(\sigma_E^{-1}\bullet\id)(u_j\tens \sigma_E^{-1}(u_i\tens e))\cr
 &=& (\alpha\,\wedge\tens\id_E\tens\id^{\tens 2})(\id\tens\sigma_E\tens\id^{\tens 2})(\sigma_E\tens\id^{\tens 3})(e\tens\xi_i\tens\xi_j\tens u_j\tens u_i)\,+ \cr
 && 
(\alpha\,\wedge\tens\id_E\tens\id)(\id\tens\sigma_E\tens\id) (\id\tens\id_E\tens\square)(\xi_i\tens \sigma_E^{-1}(u_i\tens e))  \ .
\end{eqnarray}
We shall call the two terms of the RHS of (\ref{hvcfhjjhfxxgh}) terms 4a and 4b respectively of 
(\ref{jhcfshjtfdj}). Then combining terms 1,3 and 4b gives
\begin{eqnarray}
&& (\alpha\tens\id_E\tens\id)\big(\extd\tens\id_E\tens\id-\id\wedge\nabla_E\tens\id-\id\wedge(\sigma_E\tens\id)(\id_E\tens\square)\big)  \cr
&& (\xi_i\tens\sigma_E^{-1}(u_i\tens e)) \cr
&=&  (\alpha\tens\id_E\tens\id)\big(\extd\tens\id_E\tens\id-\id\wedge\nabla_E\tens\id-\id\wedge(\sigma_E\tens\id)(\id_E\tens\square)\big) \cr
&&  (\sigma_E(e\tens \xi_i)\tens u_i)\ ,
\end{eqnarray}
where we have used the fact that the long bracket gives a well defined operator on $\Omega^1 A\tens_A E\tens_A\mathrm{Vec}A$. Now we can rewrite (\ref{hvcfhjjhfxxgh}) as
\begin{eqnarray} \label{hjgfchxfhdz}
&& \vartheta_E((\alpha\tens\id)\mathcal{R}\tens e) - (\alpha\tens\id)\mathcal{R}\,\la\,e\tens 1_A\cr
&=& -\, (\alpha\tens\id_E\tens\id)\Big((\id\wedge\sigma_E\tens\id^{\tens 2})(\sigma_E\tens\id^{\tens 3})(e\tens\xi_i\tens\xi_j\tens u_j\tens u_i) \cr
&&  +\, \Big((\extd\tens\id_E-\id\wedge\nabla_E)\,\sigma_E\tens\id - (\id\wedge\sigma_E)(\nabla_E\tens\id)\tens\id \, - \cr
&& \big(\id\wedge(\sigma_E\tens\id)\big)(\sigma_E\tens\square)
\Big)(e\tens \xi_i\tens u_i)
\Big)
\end{eqnarray}
Now we use the following equality, which is proved in \cite{bbsheaf} for vanishing curvature, but in fact the proof only requires that the curvature $R_E$ is a right $A$-module map.
\begin{eqnarray*}
\big(\extd\tens\id_E-(\id\wedge\nabla_E)\big)\sigma_E\,=\, (\id\wedge\sigma_E)(\nabla_E\tens\id)+\sigma_E(\id_E\tens\extd) \, :\, E\tens_A\Omega^1 A\to \Omega^2 A \tens_A E\ .
\end{eqnarray*}
We can use this to rewrite (\ref{hjgfchxfhdz}) as
\begin{eqnarray} \label{hjgxfhdz}
&& \vartheta_E((\alpha\tens\id)\mathcal{R}\tens e) - (\alpha\tens\id)\mathcal{R}\,\la\,e\tens 1_A\cr
&=& -\, (\alpha\tens\id_E\tens\id)\Big((\id\wedge\sigma_E\tens\id^{\tens 2})(\sigma_E\tens\id^{\tens 3})(e\tens\xi_i\tens\xi_j\tens u_j\tens u_i) \cr
&&  +\, \Big(\sigma_E(\id_E\tens\extd)\tens\id -  \big(\id\wedge(\sigma_E\tens\id)\big)(\sigma_E\tens\square)
\Big)(e\tens \xi_i\tens u_i)
\Big) \cr
&=& ((\alpha\tens\id_E)\sigma_E\tens\id)(e\tens\extd\xi_i\tens u_i - e\tens\xi_i\wedge\xi_j\tens u_j\bullet u_i) \ .
\end{eqnarray}
Now using (\ref{ajhscfxhgfz1}) gives the answer.\quad $\blacksquare$

\begin{cor}  \label{bvsdvvchas}
Suppose that $(E,\nabla_E,\sigma_E)\in {}_A\mathcal{F}_A$ , and 
 that $\Omega^2 A$ is finitely generated projective as a right $A$-module. Then
\begin{eqnarray*}
\vartheta_E(\mathcal{W}\tens_A E)\,\subset\,E\tens_A \mathcal{W}\ ,
\end{eqnarray*}
and as a result we have a well defined quotient $\vartheta_E:\mathcal{D}A\tens_A E\to
E\tens_A \mathcal{D}A$. 
\end{cor}
\noindent {\bf Proof:}\quad From Proposition \ref{bvchas}, since 
$\widehat{\mathcal{R}}\,\la\, \id_E=0$, we have 
\begin{eqnarray*}
\vartheta_E\,(\widehat{\mathcal{R}}\tens \id_E) &=& (\ev\tens\id_E\tens \id)
(\id\tens\sigma_E\tens\widehat{\mathcal{R}})(\id\tens\id_E\tens\coev_2) \cr
&& :\, (\Omega^2 A)'\tens_A E\to E\tens_A \mathcal{T}\mathrm{Vec}A\ .
\end{eqnarray*}
To get the result for the ideal under the $\bullet$ product, we use the multiplicative property of 
$\vartheta_E$ listed in (\ref{jkhzdcvjxty}).\quad $\blacksquare$

\section{Differential operators with module endomorphisms}  \label{bcadhklskuc}
We would like to consider vector fields acting on left modules $E$ with left covariant derivative $\nabla_E$, together with left module endomorphisms (in ${}_A\mathrm{End}(E)$) acting on $E$. For example, the classical Dirac operator acting on sections of the spinor bundle requires such a module map. 
 There are at least two reasons why this is not really enough. The first is that we would prefer to be able to reorder vector fields and left module maps, rather than leave them in some sort of free product. To see what happens when we do this, we use the following lemma, where $\circ$ stands for composition of functions:

\begin{lemma}
Given $(E,\nabla_E)\in{}_A\mathcal{E}$ and $T\in {}_A\mathrm{End}(E)$, we have a  left module map $\nabla(T):E\to\Omega^1A\tens_A E$ defined by
$\nabla(T)=\nabla\circ T-(\id\tens T)\circ\nabla$. 
\end{lemma}
\noindent {\bf Proof:}\quad Consider $\nabla(a.T)$ for $a\in A$.\quad$\blacksquare$

\medskip Now we have, for $T\in {}_A\mathrm{End}(E)$ and $u\in\mathrm{Vec}A$,
\begin{eqnarray} \label{bcudhskav}
\big((u\,\la)\circ T-T\circ(u\,\la)\big)(e)\ =\ (\ev\tens\id_E)(u\tens\nabla(T)(e))\ .
\end{eqnarray}
Starting from a left module map from $E$ to $E$,
we end up introducing a module map from $E$ to $\Omega^1A\tens_A E$. But there is another reason why we may wish to consider such maps.  If we were to start only with the module $E$, it might be considered as rather artificial to impose one particular covariant derivative on it. The difference between any two covariant derivatives is a left module map from $E$ to $\Omega^1A\tens_A E$. If we add such module maps from the beginning, we obtain actions on $E$ which are independent of the choice of covariant derivative on $E$. For these reasons, we make the following definition. For a left module map $S:E\to \Omega^{\tens n}A\tens_A E$ and $\underline v\in\mathrm{Vec}^{\tens n}A$ define $K_n(\underline v,S):E\to E$ by
\begin{eqnarray}
K_n(\underline v,S)(e)\ =\ (\ev^{\<n\>}\tens\id_E)(\underline v\tens S(e))\ .
\end{eqnarray}
Now (\ref{bcudhskav}) reads $(u\,\la)\circ T-T\circ(u\,\la)=K_1(u,\nabla(T))$. For completeness we define $K_0(a,T)=a.T$. To see how these operations combine, we define, for $S:E\to \Omega^{\tens n}A\tens_A E$ and $U:E\to \Omega^{\tens m}A\tens_A E$,
\begin{eqnarray} 
S\circ U\ =\ (\id^{\tens m}\tens S)U:E\to \Omega^{\tens m+n}A\tens_A E\ .
\end{eqnarray}
Now we have
\begin{eqnarray}
K_n(\underline v,S)\circ K_m(\underline w,U)\ =\ K_{n+m}(\underline v\tens \underline w,
S\circ U)\ .
\end{eqnarray}
To look at the other compositions, we will need to observe that, using (\ref{bvhkvjhgc}),
\begin{eqnarray} \label{cbdhsu}
\nabla_E(S) &=& (\square^{\<n\>}\tens\id_E+\id^{\tens n}\tens\nabla_E)S-(\sigma^{-\<n\>}\tens\id_E)(\id\tens S)\nabla_E\ .
\end{eqnarray}
is a left module map $\nabla_E(S):E\to \Omega^{\tens n+1}A\tens_A E$. To complete being able to reorder the $K_n$ to the left and the $\la$ to the right we need the following result:

\begin{propos}
For $u\in \mathrm{Vec}A$, $\underline v\in\mathrm{Vec}^{\tens n}A$
and $S:E\to \Omega^{\tens n}A\tens_A E$,
\begin{eqnarray*}
(u\,\la)\circ K_n(\underline v,S) &=& K_n(u\,\la\,\underline v,S)
+K_{n+1}(u\tens \underline v,\nabla_E(S))+ K_n(\underline v',S)\circ(u'\,\la)\ ,
\end{eqnarray*}
where $\sigma^{-\<n\>}(u\tens \underline v)= \underline v'\tens u'$ (a finite sum) and
$\sigma^{-\<n\>}:\mathrm{Vec}^{\tens n+1}A\to \mathrm{Vec}^{\tens n+1}A$ is
given by the same symbolic formula as (\ref{bvhkvjhgc}), though the domain is different. 
\end{propos}
\noindent {\bf Proof:}\quad By differentiating the evaluation and using (\ref{cbdhsu}) we have
\begin{eqnarray*}
(u\,\la)\circ K_n(\underline v,S)(e) &=& K_n(u\,\la\,\underline v,S)  \cr
&& +\ (\ev^{\<n+1\>}\tens\id_E)((u\tens \underline v)\tens
(\square^{\<n\>}\tens\id_E+\id^{\tens n}\tens\nabla_E)S(e)) \cr
&= &  K_n(u\,\la\,\underline v,S) +K_{n+1}(u\tens \underline v,\nabla_E(S))   \cr
&& +\ (\ev^{\<n+1\>}\tens\id_E)((u\tens \underline v)\tens
(\sigma^{-\<n\>}\tens\id_E)(\id\tens S)\nabla_E(e)) \ . \quad\blacksquare
\end{eqnarray*}

\section{Symbols of differential operators}  \label{bcdhlsku}
The reader will recall that classically the symbol is made from the highest order part of the differential operator. The partial derivatives are replaced by the corresponding vectors, and we get a section of some tensor power of the tangent space. Of course, we could take all orders instead of just the highest order, but then we would find that the lower order terms were coordinate dependent. The symbols have useful properties: When we compose differential operators, we take the tensor product of their symbols. Some classes of differential operators, most famously elliptic operators, have much of their behaviour determined by their symbols. 

It will now not surprise the reader that the basic idea of symbols and the composition rule work in the noncommutative setting. After all, we have been describing differential operators by taking tensor powers of vector spaces to all orders. Merely restricting to the highest order non-zero term gives the symbol. A look at Lemma~\ref{hjsvhjkvjhk} shows that composition of differential operators using the $\bullet$ product is simply $\tens_A$ on the highest order parts. Of course the lower order terms of the $\bullet$ product depend on the covariant derivative $\square$ -- changing this covariant derivative is effectively changing coordinates. 

However we need to be careful - the symbol is classically an element of a \textit{symmetric} tensor product, or equivalently expressed as a (commutative) polynomial. What we have been saying in a noncommutative context refers to $\mathcal{T}\mathrm{Vec}A_\bullet$ -- the full tensor product which acts on the category ${}_A\mathcal{E}$. In Section~\ref{bchksvkuv} 
we discuss the algebra $\mathcal{D}A$ which has extra relations imposed, and acts on  modules with zero curvature covariant derivative ${}_A\mathcal{F}$. Now these relations, taking only the highest tensor order part, apply to the symbols. For example, for the 3D calculus on deformed $SU(2)$ discussed in Section~\ref{bhsoabvi}, from (\ref{kuycvghjkjhc}) we get following relations
\begin{eqnarray} \label{kuhjkjhdscc}
 u_-\tens u_+-q^2\,u_+\tens u_-\,=\,0\ ,\quad
 u_0\tens u_\pm - q^{\pm4}\,u_\pm\tens u_0\,=\,0\ .
\end{eqnarray}
Then a symbol for $\mathcal{D}$ for this calculus really is a polynomial in $u_0,u_+,u_-$, with the relatively minimal change in product (\ref{kuhjkjhdscc}) from the classical case (taking coefficients in the algebra, if the differential operator is not left invariant). However it may be that we were lucky in this case -- more generally it might be possible that the relations on $\mathcal{D}A$ might have combinations cancelling in order 2, producing non-obvious relations in order 1. Caution is advised until the situation is understood better. 

\section{Noncommutative complex structures}  \label{bcudisokocfgy}
We shall use the version of noncommutative complex differential geometry introduced in 
\cite{BegSmiComplex} and referenced in \cite{KhLavS}. Suppose that $A$ is a star algebra with differential calculus for which the map $a.\extd b\mapsto \extd b^*.a^*$ gives a well defined star operation on $\Omega^1 A$, extending to all $\Omega^n A$. Then an almost complex structure on $A$ is a bimodule map $J:\Omega^1 M\to \Omega^1 M$ for which $J(\xi^*)=J(\xi)^*$, $J\circ J$ is minus the identity, and the map
\begin{eqnarray}
J\tens\id^{\tens n-1}+\dots+\id^{\tens n-1}\tens J:\Omega^{\tens n}A\to \Omega^{\tens n}A
\end{eqnarray}
gives a well defined map (via the $\wedge$ product) $J:\Omega^nA\to \Omega^nA$. From this we can decompose each $\Omega^nA$ into a direct sum of $\Omega^{p,q}A$, the $\mathrm{i}\,(p-q)$ eigenspace of $J$, where $p+q=n$ and $p,q\ge 0$. Take $\pi^{p,q}:\Omega^{n}A \to \Omega^{p,q}A$ to be the corresponding projection. 
The almost complex structure is called integrable if (amongst several equivalent conditions), for all 
$\xi\in \Omega^{0,1}A$ we have $\extd\xi\in \Omega^{0,2}A\oplus \Omega^{1,1}A$. 

The idea of complex structure given here is exactly the same as that used in commutative geometry -- though the star condition is ommitted as in an \textit{a priori} real manifold, it is not needed. Also the integrability condition is more often given in terms of vector fields, and we will come to this below. Now the content of the Newlander-Nirenberg integrability theorem \cite{NewNir} is that on a  real manifold, with an almost complex structure satisfying the integrability condition, we can construct local complex coordinates, linked by complex analytic transition functions. 

However we do not really know what a noncommutative `real manifold' is, and (at least with our current degree of understanding) we cannot use the idea of local coordinates. It is tempting to believe that there would be an analogue of the Newlander-Nirenberg integrability theorem in the noncommutative case, but what would replace local coordinates? In this paper I shall make an argument that we should look at covariant derivatives on the 1-forms with certain properties as a substitute for local complex coordinates. 

To deal with a few more tecnicalities and matters of notation, if $\Omega^{n}A$ is a finitely generated projective module, then so are the $\Omega^{p,q}A$ for $p+q=n$, as they sum to 
$\Omega^{n}A$. We define $\partial:\Omega^{p,q}A\to \Omega^{p+1,q}A$ and 
$\bar\partial:\Omega^{p,q}A\to \Omega^{p,q+1}A$ by $\extd$ followed by the appropriate projection, 
$\pi^{p+1,q}$ or $\pi^{p,q+1}$. By the integrability condition, $\partial^2=0$ and 
$\bar\partial^2=0$. 

As the vector fields are dual to the 1-forms (i.e.\ $\mathrm{Vec}A=\Hom_A(\Omega^1 A,A)$), many things that can be done with one can be done with the other. In particular I mentioned that the original approach to the Newlander-Nirenberg integrability theorem was via vector fields. Just to show that it can be done in a noncommutative setting, we shall now show that the integrability condition can be phrased in terms of vector fields. Remember that Definition \ref{deflie} gives our form of the Lie bracket of vector fields:

\begin{defin} \label{jkzdxcfghcjygfdz}
For an almost complex structure $J:\Omega^1 A\to \Omega^1 A$, define 
$J:\mathrm{Vec}A\to \mathrm{Vec}A$ by $J(v)=v\circ J$. Then the $\pm\mathrm{i}$ eigenspaces 
of $J$ are respectively denoted $\mathrm{Vec}^{1,0}A$ and $\mathrm{Vec}^{0,1}A$. 
\end{defin}

\begin{propos}\textbf{The Newlander-Nirenberg integrability condition.}
Suppose that the integrability condition holds for an almost complex structure on $\Omega^1 A$. 
Also suppose that $x\in \mathrm{Vec}^{1,0}A\tens \mathrm{Vec}^{1,0}A$ obeys the condition that $\ev^{\<2\>}(x\tens k)=0$ for all $k\in\ker\wedge:\Omega^{1,0}A\tens_A \Omega^{1,0}A\to \Omega^{2,0}A$. Then $\varphi(x)\in  \mathrm{Vec}^{1,0}A$.
\end{propos}
\noindent {\bf Proof:}\quad  Write $x=u\tens v$ (summation implicit). We ned to show that
$\varphi(u\tens v)(\xi)=0$ for all $\xi\in \Omega^{0,1}A$. The formula is, where $z\in \Omega^{1}A\tens_A \Omega^{1}A$ is chosen so that $\wedge z=\extd\xi$,
\[
\varphi(X\tens Y)(\xi)\,=\,D_u(v(\xi))\,+\,\ev^{\<2\>} (u\tens
v\tens z)\ ,
\]
Now $Y(\xi)=0$, so we only have to show that $\ev^{\<2\>} (u\tens
v\tens z)=0$. As $\xi\in \Omega^{0,1}A$, the differential form version of the integrability condition states that $\extd\xi\in \Omega^{0,2}A\oplus \Omega^{1,1}A$. This means that we can choose not to have any component of $z$ in $\Omega^{1,0}A\tens_A \Omega^{1,0}A$, so $\ev^{\<2\>} (u\tens v\tens z)=0$. \quad $\blacksquare$

\section{What is a holomorphic vector field?}  \label{ioygfcityk}
Classically, if we differentiate a holomorphic function along a 
holomorphic vector field, we should get another holomorphic function. In terms of local holomorphic coordinates $z_1,\dots,z_n$, the following is a holomorphic vector field, where 
the functions $f_i(z_1,\dots,z_n)$ are holomorphic functions,
\begin{eqnarray}  \label{cfjksjtdf}
v\ =\ f_1(z_1,\dots,z_n)\,\frac{\partial}{\partial z_1}+\dots+f_n(z_1,\dots,z_n)\,\frac{\partial}{\partial z_n}\ .
\end{eqnarray}
However in noncommutative geometry we do not have the luxury of working in local coordinates, we need a global definition. 

A glance at (\ref{cfjksjtdf}) will show that $v\in \mathrm{Vec}^{1,0}A$
(see Definition \ref{jkzdxcfghcjygfdz}). Now
suppose that $v\in\mathrm{Vec}^{1,0}A$ and $a\in A$ is holomorphic, i.e.\ $\pdol\,a=0$. The differential of $a$ along (or in the direction of) $v$ is $\ev(v\tens\extd a)$. For this to be holomorphic, we need $\bar\partial\, \ev(v\tens\extd a)=0$. In terms of the projection $\pi^{0,1}$ we have
$\pi^{0,1}\,\extd\, \ev(v\tens\extd a)=0$. (We remind the reader that $\pi^{p,q}$ is the projection from $\Omega^{p+q}A$ to $\Omega^{p,q}A$.) Classically, we perform this differentiation in the given coordinates, and see that the functions $f_1,\dots,f_n$ must be holomorphic. Next we say that, since the change of coordinate functions are holomorphic, it makes sense to say that the corresponding functions are holomorphic in every coordinate chart. Of course, we can't do any of this. We have to differentiate the vector field with the only tool we have for doing this, a covariant derivative. When we do this, we not surprisingly discover that there are additional conditions that we must impose on the covariant derivative, and these conditions stand in for the classical assumption that the change of coordinate functions are holomorphic. 

 To differentiate $\ev(v\tens\extd a)$ use the covariant derivatives $\square$ from (\ref{ghjksvdgwv}), and write
\begin{eqnarray}  \label{jhsagcfctyj}
\pi^{0,1}\,\extd\, \ev(v\tens\extd a) &=& (\pi^{0,1}\tens\ev)(\square\,v\tens \extd a)
+ (\ev\tens \pi^{0,1})(v\tens\square\,\extd a)\ .
\end{eqnarray}
Before continuing, it would be wise to review the meaning of (\ref{jhsagcfctyj}). In fact, it may contain no information at all. On a compact complex manifold, all holomorphic functions are constants, and (\ref{jhsagcfctyj}) becomes $0=0+0$, independently of $v$. What we are going to do is to obtain a condition which is more local in character, and which would imply $\ev(v\tens\extd a)$ holomorphic for any holomorphic $a$.

\begin{lemma}    \label{kuyfdtyytr44s}
Suppose that 
\newline a)\quad  $\wedge:\Omega^{1,0}A\tens_A \Omega^{0,1}A\to \Omega^{1,1}A$
is an isomorphism with inverse $\phi$,
\newline b) \quad $(J\tens\id)\square=\square\,J:\Omega^1 A\to \Omega^{1}A\tens_A \Omega^{1}A$,
\newline c) \quad $\pi^{1,1}\mathrm{Tor}_R:\Omega^1 A\to \Omega^{1,1}A$
vanishes.
\newline Then we have, for $v\in\mathrm{Vec}^{1,0}A$ and $a\in A$,
\begin{eqnarray*} 
\pi^{0,1}\,\extd\, \ev(v\tens\extd a) 
&=& (\pi^{0,1}\tens\ev)(\square\,v\tens \partial a)
- (\ev\tens \id)(v\tens \phi\,\pi^{1,1}\wedge\square\,\bar\partial a)\ .
\end{eqnarray*}
\end{lemma}
\noindent {\bf Proof:}\quad To repeat (\ref{jhsagcfctyj}), for $v\in\mathrm{Vec}^{1,0}A$
and $a\in A$,
\begin{eqnarray*} 
\pi^{0,1}\,\extd\, \ev(v\tens\extd a) &=&\pi^{0,1}\,\extd\, \ev(v\tens\partial a) \cr
&=& (\pi^{0,1}\tens\ev)(\square\,v\tens \partial a)
+ (\ev\tens \pi^{0,1})(v\tens\square\,\partial a)\ .
\end{eqnarray*}
By \textit{(b)}
\begin{eqnarray*}
\square\,\partial a &=& (\id\tens \pi^{0,1})\,\square\,\partial a + 
(\id\tens \pi^{1,0})\,\square\,\partial a\ \in\, \Omega^{1,0}A\tens_A \Omega^{0,1}A \bigoplus \Omega^{1,0}A\tens_A \Omega^{1,0}A\ .
\end{eqnarray*}
From this it follows that, where $\phi$ is the inverse of
 $\wedge:\Omega^{1,0}A\tens_A \Omega^{0,1}A\to \Omega^{1,1}A$. 
\begin{eqnarray*}
\pi^{1,1}\wedge\square\,\partial a &=& \wedge(\id\tens \pi^{0,1})\,\square\,\partial a \ ,\cr
\phi\,\pi^{1,1}\wedge\square\,\partial a &=& (\id\tens \pi^{0,1})\,\square\,\partial a \ ,
\end{eqnarray*}
Now remember the definition (\ref{tordefee}) of $\mathrm{Tor}_R$, and using the fact that $\extd^2a=0$ we get
\begin{eqnarray*}
0\ =\ \pi^{1,1}\mathrm{Tor}_R(\extd a) &=& \pi^{1,1}\wedge\,\square\,\extd a\ =\ 
\pi^{1,1}\wedge\,\square\,\partial a+ \pi^{1,1}\wedge\,\square\,\bar\partial a \ ,
\end{eqnarray*}
giving $(\id\tens \pi^{0,1})\,\square\,\partial a=-\phi\,\pi^{1,1}\wedge\square\,\bar\partial a$
and the result.\quad$\blacksquare$

\begin{cor}  \label{kuyfdtyytrvsv44s} Given the conditions of Lemma \ref{kuyfdtyytr44s},
 if $(\pi^{0,1}\tens\id)\square\,v=0$ for $v\in\mathrm{Vec}^{1,0}A$, it follows that if $a\in A$ is holomorphic, then so is $\ev(v\tens\extd a)$. 
\end{cor}

Returning to the conditions of Proposition \ref{kuyfdtyytr44s} from the point of view of examples, it is instructive to note the assumption that $\wedge:\Omega^{p,0}A\tens_A \Omega^{0,q}A\to \Omega^{p,q}A$ is an isomorphism is not only true for the noncommutative complex projective spaces in 
\cite{buncomproj}, but that the assumption is fundamental to the whole construction there. 
Of course, classically condition \textit{(a)}  is just reordering wedge products of coordinate functions using $\extd\bar z_i\wedge\extd z_j=-\extd z_j\wedge\extd\bar z_i$. 

\section{Nice covariant derivatives on integrable complex structures} \label{bcdsuovuy}
To review Section~\ref{ioygfcityk}, we should ask why the conditions on $\square$ in Lemma \ref{kuyfdtyytr44s} are required. It should be remembered that the form of the holomorphic vector field in (\ref{cfjksjtdf}) 
and its behaviour is completely determined by our idea of local complex coordinates in the classical case. In the noncommutative setting we lose this, and all we have to replace it is $\square$. In other words, $\square$ needs to contain all the information about the local complex coordinates, and their analytic transition functions, that we need to complete the proof. From this point of view, it is not surprising that $\square$ will be subject to various conditions, as it is basically encoding much of the idea of a `noncommutative complex manifold'. To continue the discussion from Section~\ref{bcudisokocfgy}, the conditions of Lemma \ref{kuyfdtyytr44s} might form part of what replaces `local holomorphic coordinates' in noncommutative geometry. It is possible to go part way to building a `nice' covariant derivative that has these properties, as follows. 

To begin, remember that any finitely generated projective right $A$-module $E$ can be given a right covariant derivative $\tilde\nabla:E\to E\tens_A\Omega^1 A$. Given a dual basis $e_i\in E$ and $e^i\in \Hom_A(E,A)$, define $\tilde\nabla(e_i)=e_j\tens\Gamma_{ji}$ where we consider $\Gamma_{ji}$ as a matrix with entries in $\Omega^1 A$. The relations between the $e_i$ are given by $e_i\,P_{ij}=e_j$, where $P_{ij}=e^i(e_j)$ is a projection matrix with values in $A$. Without loss of generality we may assume that $P\,\Gamma=\Gamma$ (using the matrix product). 
Applying $\tilde\nabla$ to the relation gives the required equation $P\,\Gamma\,(1-P)=P.\extd P$. This has solutions $\Gamma=P.\extd P+PBP$, where $B$ is any matrix with entries in $\Omega^1 A$. Of course, constructing a right bimodule covariant derivative is more complicated, if it is possible at all. 

Now $\Omega^{1,0}A$ and $\Omega^{0,1}A$ are direct summands of $\Omega^1 A$, so if $\Omega^1 A$ is right finitely generated projective it follows that both $\Omega^{1,0}A$ and $\Omega^{0,1}A$ are also right finitely generated projective. We shall suppose that we are given right bimodule covariant derivatives $\nabla_{10}:\Omega^{1,0}A\to \Omega^{1,0}A\tens_A \Omega^{1}A$ and $\nabla_{01}:\Omega^{0,1}A\to \Omega^{0,1}A\tens_A \Omega^{1}A$, with corresponding $\sigma_{10}^{-1}$ and $\sigma_{01}^{-1}$. It might be possible to further assume that $\nabla_{10}$ and $\nabla_{01}$ are linked by conjugation $\star:\Omega^{1,0}A\to\overline{\Omega^{0,1}A}$, but we shall not be concerned by this now. What does concern us is that we can construct a right bimodule covariant derivative on $\Omega^1A$ just by adding the ones on $\Omega^{1,0}A$ and $\Omega^{0,1}A$. However we shall modify 
$\nabla_{10}$ and $\nabla_{01}$ first. 

\begin{propos}  \label{cbyiaokivy}
Suppose that we have right bimodule covariant derivatives $(\Omega^{1,0}A,\nabla_{10},\sigma_{10}^{-1})$ and 
$(\Omega^{0,1}A,\nabla_{01},\sigma_{01}^{-1})$. We suppose that
\newline a)\quad  $\wedge:\Omega^{1,0}A\tens_A \Omega^{0,1}A\to \Omega^{1,1}A$
is an isomorphism with inverse $\phi$,
\newline b)\quad  $\wedge:\Omega^{0,1}A\tens_A \Omega^{1,0}A\to \Omega^{1,1}A$
is an isomorphism with inverse $\psi$. 
\newline Then there is a right bimodule covariant derivative $(\Omega^1 A,\square,\sigma^{-1})$ defined by
\begin{eqnarray*}
\square &=& \big((\id\tens\pi^{1,0})\nabla_{10}-\phi\,\pi^{1,1}\extd\big)\,\pi^{1,0}
+ \big((\id\tens\pi^{0,1})\nabla_{01}-\psi\,\,\pi^{1,1}\extd\big)\,\pi^{0,1}\ ,\cr
\sigma^{-1} &=& \big((\id\tens\pi^{1,0})\sigma^{-1}_{10}-\phi\,\pi^{1,1}\wedge\big)(\id\tens\pi^{1,0})
+ \big((\id\tens\pi^{0,1})\sigma^{-1}_{01}-\psi\,\pi^{1,1}\wedge\big)(\id\tens\pi^{0,1})\ .
\end{eqnarray*}

Further $(\Omega^1 A,\square,\sigma^{-1})$ satisfies the properties:
\newline i) \quad $(J\tens\id)\square=\square\,J:\Omega^1 A\to \Omega^{1}A\tens_A \Omega^{1}A$,
\newline ii) \quad $\pi^{1,1}\mathrm{Tor}_R:\Omega^1 A\to \Omega^{1,1}A$ vanishes.
\newline iii)  \quad $
\pi^{1,1}\wedge\sigma^{-1} = -\ \pi^{0,1}\wedge\pi^{1,0}-\pi^{1,0}\wedge \pi^{0,1}
:\Omega^{\tens 2}A\to \Omega^{1,1}A$\ .

\end{propos}
\noindent {\bf Proof:}\quad First check the right Leibniz rule, for $\xi\in\Omega^{1,0}A$ and $\eta\in\Omega^{0,1}A$,
\begin{eqnarray*}
\square(\xi.a)-\square(\xi).a &=& \xi\tens\partial a + \phi\,\pi^{1,1}(\xi\wedge\extd a)\ =\ 
 \xi\tens\partial a +  \xi\tens\bar\partial a\,=\,\xi\tens\extd a \ , \cr
 \square(\eta.a)-\square(\eta).a &=& \eta\tens\bar\partial a + \psi\,\pi^{1,1}(\eta\wedge\extd a)\ =\ 
 \eta\tens\bar\partial a + \eta\tens\partial a\ =\ \eta\tens\extd a   \ .
\end{eqnarray*}
Now check the bimodule covariant derivative rule:
\begin{eqnarray*}
\square(a.\xi)-a.\square(\xi) &=& (\id\tens\pi^{1,0})\sigma^{-1}_{10}(\extd a\tens\xi)-\phi\,\pi^{1,1}(\extd a\wedge\xi) \ , \cr
 \square(a.\eta)-a.\square(\eta) &=& (\id\tens\pi^{0,1})\sigma^{-1}_{01}(\extd a\tens\xi)-\psi\,\pi^{1,1}(\extd a\wedge\eta)   \ .
\end{eqnarray*}
Next \textit{(i)} works by looking at $\Omega^{1,0}A$ and $\Omega^{0,1}A$, which are eigenspaces for $J$. For 
\textit{(ii)},
\begin{eqnarray*}
\pi^{1,1}\wedge\square(\xi) &=& \pi^{1,1}(\id\wedge\pi^{1,0})\nabla_{10}(\xi)-\pi^{1,1}\extd\xi\ =\ -\,\pi^{1,1}\extd\xi\ ,\cr
\pi^{1,1}\wedge\square(\eta) &=& \pi^{1,1}(\id\wedge\pi^{0,1})\nabla_{01}(\eta)-\pi^{1,1}\extd\eta\ =\ -\,\pi^{1,1}\extd\eta\ .
\end{eqnarray*}
For \textit{(iii)}, we get
\begin{eqnarray*}
\pi^{1,1}\wedge\sigma^{-1} &=& (-\pi^{1,1}\wedge)(\id\tens\pi^{1,0})
+ (-\pi^{1,1}\wedge)(\id\tens\pi^{0,1}) \cr
&=& -\ \pi^{0,1}\wedge\pi^{1,0}-\pi^{1,0}\wedge \pi^{0,1}\ .\quad\blacksquare
\end{eqnarray*}

\medskip It is likely that examples could be constructed from complex structures on the noncommutative torus \cite{PS1}, but here is one on the Podle\'s sphere.

\begin{example}
There is a Hopf algebra map $\pi$ from quantum $SL_2$ (see Section~\ref{bhsoabvi}) to 
$\mathbb{CZ}$ (the group
algebra of the group $(\mathbb{Z},+)$, with group 
generator $z$) given by $\pi(a)=z$,  $\pi(b)=\pi(c)=0$ and $\pi(d)=z^{-1}$. This is used to construct a right coaction $(\id\tens\pi)\Delta$ of $\mathbb{CZ}$ on quantum $SL_2$, and the standard Podle\'s sphere $S^2_q$ \cite{Pod} is given as the invariants of this coaction. This coaction extends to the 3D calculus on quantum $SL_2$, by $e^\pm\mapsto e^\pm\tens z^{\pm2}$ and $e^0\mapsto e^0\tens 1$. The horizontal invariant forms on quantum $SL_2$ under this coaction give a differential calculus on $S^2_q$. The 1-forms are then elements of 
quantum $SL_2$ times $e^\pm$ which are invariant to the right $\mathbb{CZ}$ coaction. The differential calculus is given by setting $\extd e^\pm=0$, and there is a trivial $\Omega^2S^2_q$ with generator $e^+\wedge e^-$. In \cite{Maj} the spin geometry of $S^2_q$ was studied,
splitting $\Omega^1S^2_q$ into its $e^+$ and $e^-$ components.

An integrable almost complex structure is given by $J(e^\pm)=\pm\,\mathrm{i}\,e^\pm$. We take a right covariant derivative $\square$ on $\Omega^1 S^2_q$ defined by $\square(e^\pm)=0$. This is torsion free, and satisfies the conditions of Lemma~\ref{kuyfdtyytr44s}. 
We can also take projections to appropriate eigenspaces of $J$ to define $\nabla_{10}$ and 
$\nabla_{01}$ satisfying the conditions of Proposition~\ref{cbyiaokivy}. 
\end{example}

\section{The subalgebras $\mathcal{T}\mathrm{Vec}^{*,0}A_\bullet$ and $\mathcal{T}\mathrm{Vec}^{0,*}A_\bullet$} \label{vukjhcfxfth}
We consider the subspaces of $\mathcal{T}\mathrm{Vec}A$,
\begin{eqnarray*}
\mathcal{T}\mathrm{Vec}^{*,0}A &=& A \oplus \mathrm{Vec}^{1,0}A\oplus (\mathrm{Vec}^{1,0}A\tens_A \mathrm{Vec}^{1,0}A)\oplus\dots \cr
\mathcal{T}\mathrm{Vec}^{0,*}A &=& A \oplus \mathrm{Vec}^{0,1}A\oplus (\mathrm{Vec}^{0,1}A\tens_A \mathrm{Vec}^{0,1}A)\oplus\dots
\end{eqnarray*}
and ask when they are subalgebras of $\mathcal{T}\mathrm{Vec}A_\bullet$ under the $\bullet$ product. We use the notation $\Omega^{\tens n,0}A=(\Omega^{1,0}A)^{\tens n}$,
$\Omega^{\tens 0,n}A=(\Omega^{0,1}A)^{\tens n}$,
 $\mathrm{Vec}^{\tens n,0}A=(\mathrm{Vec}^{1,0}A)^{\tens n}$ and $\mathrm{Vec}^{\tens 0,n}A=(\mathrm{Vec}^{0,1}A)^{\tens n}$.

\begin{lemma}  \label{bchadjksvcjhgc}
Suppose that $\square:\Omega^1 A\to \Omega^1 A\tens_A \Omega^1 A$ is a right covariant derivative with $(J\tens\id)\square=\square\,J$. Then the dual left covariant derivative $\square:\mathrm{Vec} A\to \Omega^1 A\tens_A \mathrm{Vec}A$ obeys $(\id\tens J)\square=\square\,J$. If $\square$ is also a bimodule covariant derivative, then 
\begin{eqnarray*}(J\tens \id)\,\sigma^{-1} &=&  \sigma^{-1}\,(\id\tens J): \Omega^{\tens 2} A
\to \Omega^{\tens 2} A\ ,\cr
(\id\tens J)\,\sigma &=&  \sigma\,(J\tens\id):  \mathrm{Vec}A \tens_A \Omega^1 A
\to \Omega^1 A\tens_A \mathrm{Vec}A\ .
\end{eqnarray*}
\end{lemma}
\noindent {\bf Proof:}\quad By definition of $J$ on $\mathrm{Vec}A$ we have
$\ev(J\tens\id)=\ev(\id\tens J)$ on $\mathrm{Vec}A\tens_A \Omega^1 A$. Differentiating this
using (\ref{ghjksvdgwv}) gives the first result. For the last equation (the one before is similar),
\begin{eqnarray*}
\sigma(J(v)\tens\extd a) &=& \square(J(v).a)-\square(J(v)).a \cr
&=& \square(J(v.a))-\square(J(v)).a \cr
&=& (\id\tens J)(\square(v.a)-\square(v).a)\ .\quad\blacksquare
\end{eqnarray*}

\begin{propos} \label{bchdisvoi}
Suppose that $\square:\Omega^1 A\to \Omega^1 A\tens_A \Omega^1 A$ is a right bimodule covariant derivative with $(J\tens\id)\square=\square\,J$. Then for all $m\ge 0$,
\begin{eqnarray*}
\square^{\<m\>}\underline \xi\in \Omega^{\tens m,0}A\tens_A  \Omega^1A& \forall
\ \underline \xi\in \Omega^{\tens m,0}A\ ,\cr
\square^{\<m\>}\underline \xi\in \Omega^{\tens 0,m}A\tens_A  \Omega^1A& \forall
\ \underline \xi\in \Omega^{\tens 0,m}A\ ,\cr
\square^{\<m\>}\underline w\in \Omega^1A\tens_A\mathrm{Vec}^{\tens m,0}A & \forall
\ \underline w\in \mathrm{Vec}^{\tens m,0}A\ ,\cr
\square^{\<m\>}\underline w\in\Omega^1A\tens_A \mathrm{Vec}^{\tens 0,m}A & \forall
\ \underline w\in \mathrm{Vec}^{\tens 0,m}A\ .
\end{eqnarray*}
\end{propos}
\noindent {\bf Proof:}\quad For $m=0$ the result is immediate.  The $m=1$ cases are given by the recursive definitions
(\ref{vcakkccjjcj}) and (\ref{kjvycyjjhgx}),  and Lemma~\ref{bchadjksvcjhgc}
 together with the definitions of $\Omega^{\tens 1,0}A$ etc.\ as eigenspaces of $J$. 
Now we shall prove the third equation by induction (the others are similar). Suppose that the third equation is true for $m$. 
 From (\ref{kjvycyjjhgx}) we have
\begin{eqnarray*} 
\square^{\<m+1\>} &=& \square\tens \id^{\tens m}+
(\sigma\tens\id^{\tens n})(\id\tens\square^{\<m\>})\ .
\end{eqnarray*}
Now Lemma~\ref{bchadjksvcjhgc} gives
\begin{eqnarray*}
\sigma:\mathrm{Vec}^{1,0}A\tens_A \Omega^1 A\to \Omega^1 A\tens_A \mathrm{Vec}^{1,0}A\ ,
\end{eqnarray*}
and this gives the result for $\square^{\<m+1\>}$. \quad$\blacksquare$

\begin{cor}   \label{bvuislvdiubs}
Suppose that $\square:\Omega^1 A\to \Omega^1 A\tens_A \Omega^1 A$ is a right bimodule covariant derivative with $(J\tens\id)\square=\square\,J$. Then 
$\mathcal{T}\mathrm{Vec}^{*,0}A_\bullet$ and $\mathcal{T}\mathrm{Vec}^{0,*}A_\bullet$ are subalgebras of $\mathcal{T}\mathrm{Vec}A_\bullet$ under the $\bullet$ product. 
\end{cor}
\noindent {\bf Proof:}\quad 
From Lemma \ref{hjsvhjkvjhk} (following \cite{bbdiff}) we see that for $\mathcal{T}\mathrm{Vec}^{*,0}A$ to be a subalgebra of $\mathcal{T}\mathrm{Vec}A_\bullet$,
 all we need is that, for all $u\in \mathrm{Vec}^{1,0}A$ and $\underline w\in \mathrm{Vec}^{\tens m,0}A$,
\begin{eqnarray*} 
(\ev\tens\id^{\tens m})(u\tens\square^{\<m\>}\underline w)\in \mathrm{Vec}^{\tens m,0}A\ .
\end{eqnarray*}
This is proved in Proposition~\ref{bchdisvoi}. The other case is similar. \quad$\blacksquare$

\medskip
Since $\mathcal{T}\mathrm{Vec}^{*,0}A_\bullet$ and $\mathcal{T}\mathrm{Vec}^{0,*}A_\bullet$ are subalgebras of $\mathcal{T}\mathrm{Vec}A_\bullet$, they act by the usual formula on objects in ${}_A\mathcal{E}$. However they act on objects in other categories as well:

\begin{defin}
The category ${}_A\mathcal{H}$ consists of objects $(E,\partial_E)$, where $E$ is a left $A$-module, and $\partial_E:E\to\Omega^{1,0}A\tens_A E$ obeys the $\partial$ Leibniz rule
$\partial_E(a.e)=\partial a\tens e+a.\partial_E(e)$ for all $a\in A$ and $e\in E$. The morphisms
$T:(E,\partial_E)\to (F,\partial_F)$ are left $A$-module maps $T:E\to F$ for which $\partial_F\,T=(\id\tens T)\partial_E$. 

The category ${}_A\bar{\mathcal{H}}$ consists of objects $(E,\bar\partial_E)$, where $E$ is a left $A$-module, and $\bar\partial_E:E\to\Omega^{0,1}A\tens_A E$ obeys the $\bar\partial$ Leibniz rule
$\bar\partial_E(a.e)=\bar\partial a\tens e+a.\bar\partial_E(e)$ for all $a\in A$ and $e\in E$.  The morphisms
$T:(E,\bar\partial_E)\to (F,\bar\partial_F)$ are left $A$-module maps $T:E\to F$ for which $\bar\partial_F\,T=(\id\tens T)\bar\partial_E$. 

There are `forgetful' functors $h:{}_A\mathcal{E}\to {}_A\mathcal{H}$ and
$\bar h:{}_A\mathcal{E}\to {}_A\bar{\mathcal{H}}$ given by $h(E,\nabla_E)=(E,(\pi^{1,0}\tens\id_E)\nabla_E)$ and $\bar h(E,\nabla_E)=(E,(\pi^{0,1}\tens\id_E)\nabla_E)$. 
\end{defin}

As in (\ref{jktycfxtu}) with $\nabla_E$, 
we can iterate $\partial_E$ to get $\partial_E^{(n)}:E\to \Omega^{\tens n,0}A\tens_A E$,
and $\bar\partial_E$ to get $\bar\partial_E^{(n)}:E\to \Omega^{\tens 0,n}A\tens_A E$
recursively by the following, assuming the conditions of Corollary \ref{bvuislvdiubs},
\begin{eqnarray}\label{jktycfxtuxx}
\partial_E^{(1)} &=& \partial_E\ ,\quad \partial_E^{(n+1)}=((\id^{\tens n}\tens\pi^{1,0})\square^{\<n\>}\tens\id_E+\id^{\tens n}\tens\partial_E)\,\partial_E^{(n)}\ ,\cr
\bar\partial_E^{(1)} &=& \bar\partial_E\ ,\quad \bar\partial_E^{(n+1)}=((\id^{\tens n}\tens\pi^{0,1})\square^{\<n\>}\tens\id_E+\id^{\tens n}\tens\bar\partial_E)\,\bar\partial_E^{(n)}\ .
\end{eqnarray}
Now the algebra $\mathcal{T}\mathrm{Vec}^{*,0}A_\bullet$ acts on objects in ${}_A\mathcal{H}$ by the usual looking formula $\underline v\,\la\,e=(\ev^{\<n\>}\tens\id_E)(\underline v\tens \partial_E^{(n)} e) $. To avoid confusion we note that the element $\underline v\in\mathcal{T}\mathrm{Vec}^{*,0}A_\bullet$ acts in exactly the same way on $(E,\nabla_E)\in{}_A{\mathcal{E}}$ as on $h(E,\nabla_E)\in{}_A{\mathcal{H}}$, under the conditions of  Lemma \ref{bchadjksvcjhgc}.
 Of course, the corresponding comments hold for $\mathcal{T}\mathrm{Vec}^{0,*}A_\bullet$ acting on objects in ${}_A\bar{\mathcal{H}}$. To show this we use the following result:

\begin{lemma}  \label{bchadjvdsvvcjhgc}
Suppose $\square:\Omega^1 A\to \Omega^1 A\tens_A \Omega^1 A$ is a right covariant derivative with $(J\tens\id)\square=\square\,J$. If $(E,\partial_E)=h(E,\nabla_E)$
and $(E,\bar\partial_E)=\bar h(E,\nabla_E)$ for $(E,\nabla_E)\in{}_A{\mathcal{E}}$, then
\begin{eqnarray*}
\partial_E^{(n)} &=& \big((\pi^{1,0})^{\tens n}\tens\id_E\big)\,\nabla_E^{(n)}\ ,\cr
\bar\partial_E^{(n)} &=& \big((\pi^{0,1})^{\tens n}\tens\id_E\big)\,\nabla_E^{(n)}\ .
\end{eqnarray*}
\end{lemma}
\noindent {\bf Proof:}\quad We shall only consider the first equation. 
Refer to  (\ref{jktycfxtu}) for the definition of $\nabla_E^{(n)}$. The $n=1$ case is just the definition of 
$h(E,\nabla_E)$. Now suppose that the equation is true for $n$, and consider
\begin{eqnarray*}
\partial_E^{(n+1)} &=& ((\id^{\tens n}\tens\pi^{1,0})\square^{\<n\>}\tens\id_E+\id^{\tens n}\tens\partial_E)\,\partial_E^{(n)}  \cr
 &=& \big((\id^{\tens n}\tens\pi^{1,0})\square^{\<n\>}(\pi^{1,0})^{\tens n}\tens\id_E+(\pi^{1,0})^{\tens n}\tens\partial_E\big)\,\nabla_E^{(n)}\ .
\end{eqnarray*}
From this we can complete the proof by induction, if we have the following statement:
\begin{eqnarray} \label{cbiosuyiuoyv}
\square^{\<n\>}(\pi^{1,0})^{\tens n}\ =\ ((\pi^{1,0})^{\tens n}\tens\id)\,\square^{\<n\>}
:\Omega^{\tens n}A\to \Omega^{\tens n+1}A\ .
\end{eqnarray}
We now prove (\ref{cbiosuyiuoyv}) by induction. As $\pi^{1,0}$ can be written in terms of  $J$,
Lemma~\ref{bchadjksvcjhgc} gives the $n=1$ case. Now suppose that the equation (\ref{cbiosuyiuoyv}) is true for $n$, and consider, from (\ref{vcakkccjjcj})
and Lemma~\ref{bchadjksvcjhgc},
\begin{eqnarray} \label{vcakkjcj}
\square^{\<n+1\>}(\pi^{1,0})^{\tens n+1} &=& (\pi^{1,0})^{\tens n}\tens\square\,\pi^{1,0}+(\id^{\tens n}\tens\sigma^{-1})(\square^{\<n\>}(\pi^{1,0})^{\tens n}\tens\pi^{1,0})\cr
&=& \big((\pi^{1,0})^{\tens n+1}\tens\id\big)\big(\id^{\tens n}\tens\square \big) \cr 
&& +\,\big((\pi^{1,0})^{\tens n}\tens\sigma^{-1}(\id\tens\pi^{1,0})\big)(\square^{\<n\>}\tens\id)\cr
&=& ((\pi^{1,0})^{\tens n+1}\tens\id)\,\square^{\<n+1\>} \ .\quad\blacksquare
\end{eqnarray}

\begin{example} \label{iouycgziuyf1}
As an example of an object in ${}_A\mathcal{H}$,
we can choose $\mathcal{T}\mathrm{Vec}^{*,0}A_\bullet$, with covariant derivative
\begin{eqnarray*}
\partial_{\mathcal{HD}}(\underline v)\ =\ \eta_i\tens\eta^i\bullet\underline v\ ,
\end{eqnarray*}
where we take  a dual basis $\eta_i\tens\eta^i\in\Omega^{1,0}A\tens_A \mathrm{Vec}^{1,0}A$. 
This is simply the holomorphic analogue of (\ref{kjhsacvkyavkuxc}) for $\mathcal{T}\,\mathrm{Vec} A_\bullet$. 
Also, just as for (\ref{kjhsacvkyavkuxc}), $\partial_{\mathcal{HD}}$ is a right module map, and so we have a bimodule covariant derivative in a very trivial manner. 
The reader can now construct a corresponding example of an object in ${}_A\bar{\mathcal{H}}$. 
This is analogous to the classical construction given in \cite{MasSabDmod}.
\end{example}

\begin{example} \label{iouycgziuyf2}
 \cite{BegSmiComplex} Assume that  $\wedge:\Omega^{1,0}A\tens_A \Omega^{0,1}A\to \Omega^{1,1}A$
is an isomorphism with inverse $\phi$. Then we can take $\Omega^{0,1}A$ as an object in ${}_A\mathcal{H}$ by equipping it with the left $\partial$ covariant derivative $\phi\,\pi^{1,1}\,\extd:\Omega^{0,1}A\to \Omega^{1,0}A\tens_A \Omega^{0,1}A$. This is also a bimodule covariant derivative with
$\sigma(\xi\tens\eta)=-\phi\,\pi^{1,1}\,(\xi\wedge\eta)$. 

\end{example}

\section{Holomorphic products of vector fields} \label{nbcilhadskuv}
The reader will recall how in Corollary \ref{kuyfdtyytrvsv44s} we showed that holomorphic vector fields $v\in\mathrm{Vec}^{1,0}A$ could be given by the `obvious' formula $(\pi^{0,1}\tens\id)\square v=0$, but at the cost of imposing various conditions on $\square$. Here we shall consider what it means for an element of $\mathcal{T}\mathrm{Vec}^{*,0}A_\bullet$ to be holomorphic. Our guiding principle is the same as that of Corollary~\ref{kuyfdtyytrvsv44s}:
$\underline{v}\in \mathcal{T}\mathrm{Vec}^{*,0}A_\bullet$ should be holomorphic if $\bar\partial(\underline{v}\,\la\, a)=0$ whenever $\bar\partial a=0$. 
The method is to introduce yet another covariant derivative
$\bar\partial_{\mathcal{HD}}:\mathcal{T}\mathrm{Vec}^{*,0}A_\bullet\to \Omega^{0,1}A\tens_A \mathcal{T}\mathrm{Vec}^{*,0}A_\bullet$.  For notation we will use $\mathrm{Vec}^{\tens \le n,0}A$ to denote the subspace of 
$\mathcal{T}\mathrm{Vec}^{*,0}A$ consisting of the direct sum of all $\mathrm{Vec}^{\tens m,0}A$
for $0\le m\le n$. 

\begin{propos} \label{cbhadjsvhjgcgfx}
Suppose that the right bimodule covariant derivative $(\Omega^1 A,\square,\sigma^{-1})$ satisfies the properties:
\newline a)\quad  $\wedge:\Omega^{1,0}A\tens_A \Omega^{0,1}A\to \Omega^{1,1}A$
is an isomorphism with inverse $\phi$,
\newline b) \quad $(J\tens\id)\square=\square\,J:\Omega^1 A\to \Omega^{1}A\tens_A \Omega^{1}A$,
\newline c) \quad $\pi^{1,1}\mathrm{Tor}_R:\Omega^1 A\to \Omega^{1,1}A$ vanishes.
\newline d)  \quad $
\pi^{1,1}\wedge\sigma^{-1} = -\ \pi^{0,1}\wedge\pi^{1,0}-\pi^{1,0}\wedge \pi^{0,1}$\ .

Define $\bar\partial_{\mathcal{HD}}:\mathrm{Vec}^{\tens \le n,0}A\to \Omega^{0,1}A\tens_A \mathrm{Vec}^{\tens \le n,0}A$ recursively, beginning with $\bar\partial_{\mathcal{HD}}=\bar\partial:A\to \Omega^{0,1}A$ for $n=0$ and
$\bar\partial_{\mathcal{HD}}=(\pi^{0,1}\tens\id)\square:\mathrm{Vec}^{1,0}A\to \Omega^{0,1}A\tens_A \mathrm{Vec}^{1,0}A$. Suppose
$u\in \mathrm{Vec}^{1,0}A$ and $\underline v\in \mathrm{Vec}^{\tens n,0}A$,
with $\xi\tens \underline w=\bar\partial_{\mathcal{HD}}( \underline v)$, use a dual basis $\eta_i\tens\eta^i\in\Omega^{1,0}A\tens_A \mathrm{Vec}^{1,0}A$, and set
\begin{eqnarray*}
\bar\partial_{\mathcal{HD}}( u\bullet\underline v) &=& ((\pi^{0,1}\tens\id)\square\,u)\bullet \underline v
- (\ev\tens \id)(u\tens \phi(\pi^{0,1}\wedge\pi^{1,0})\square(\xi))\tens \underline w \cr
&&  - \ (\ev\tens \id)(u\tens \phi(\xi\wedge\eta_i))\tens \eta^i\bullet \underline w\ .
\end{eqnarray*}
Then the following properties are true:
\newline\noindent i)\quad $\bar\partial_{\mathcal{HD}}:\mathrm{Vec}^{\tens \le n,0}A\to \Omega^{0,1}A\tens_A \mathrm{Vec}^{\tens \le n,0}A$ is well defined
\newline\noindent ii)\quad $\bar\partial_{\mathcal{HD}}:\mathrm{Vec}^{\tens \le n,0}A\to \Omega^{0,1}A\tens_A \mathrm{Vec}^{\tens \le n,0}A$ is a left $\bar\partial$ covariant derivative
\end{propos}
\noindent {\bf Proof:}\quad The properties \textit{(i-ii)} are true for $n=1$. 
Suppose that \textit{(i-ii)} are true for $n$. 

There are two parts for showing  \textit{(i)} for $n+1$. First we show that the displayed formula in the hypothesis gives a unique answer for 
$u\in \mathrm{Vec}^{1,0}A$ and $\underline v\in \mathrm{Vec}^{\tens n,0}A$. Next we show that this gives a well defined function on $\mathrm{Vec}^{\tens \le n+1,0}A$. 

For the first part, begin by noting that $\id_E\tens\bullet:E\tens_A \mathrm{Vec}A\tens \mathrm{Vec}^{\tens \le n,0}A\to E\tens_A \mathrm{Vec}^{\tens \le n+1,0}A$ is well defined (for any right $A$ module $E$), so the $\bullet$ in the displayed equation do not give rise to an ambiguity. The other source of ambiguity is that we only know 
$\xi\tens \underline w=\bar\partial_{\mathcal{HD}}( \underline v)\in\Omega^{0,1}A\tens_A \mathrm{Vec}^{\tens n,0}A$ (with the emphasis on the $\tens_A$). This means that if we were to add $\eta.a\tens\underline{x}-\eta\tens a.\underline{x}$ (where $a\in A$) to $\xi\tens \underline w$ that it should make no difference to the answer. Now
\begin{eqnarray}
&& \phi(\pi^{0,1}\wedge\pi^{1,0})\square(\eta.a)\tens \underline x + \phi(\eta.a\wedge\eta_i)\tens \eta^i\bullet \underline x  \cr
&=&  \phi(\pi^{0,1}\wedge\pi^{1,0})\square(\eta).a\tens \underline x +  \phi(\eta\wedge\partial a)\tens \underline x + \phi(\eta.a\wedge\eta_i)\tens \eta^i\bullet \underline x \ , \cr
&& \phi(\pi^{0,1}\wedge\pi^{1,0})\square(\eta).a\tens \underline x + \phi(\eta\wedge\eta_i)\tens \eta^i\bullet (a.\underline x)  \cr
&=& \phi(\pi^{0,1}\wedge\pi^{1,0})\square(\eta).a\tens \underline x 
 + \phi(\eta\wedge\eta_i).\eta^i(\partial a)\tens \underline x + \phi(\eta\wedge\eta_i)\tens (\eta^i.a)\bullet \underline x  \cr
 &=& \phi(\pi^{0,1}\wedge\pi^{1,0})\square(\eta).a\tens \underline x 
 + \phi(\eta\wedge \partial a).\underline x + \phi(\eta\wedge\eta_i)\tens (\eta^i.a)\bullet \underline x  \ ,
\end{eqnarray}
and now we use the fact that $\eta_i\tens\eta^i.a=a.\eta_i\tens\eta^i\in \Omega^{1,0}A\tens_A \mathrm{Vec}^{1,0}A$. 

For the second part of \textit{(i)} for $n+1$
 we must prove the equality, for $a\in A$
\begin{eqnarray}  \label{bhadksgc}
\bar\partial_{\mathcal{HD}}(u\bullet(a.\underline v)) &=& \bar\partial_{\mathcal{HD}}((u.a)\bullet \underline v) +
 \bar\partial_{\mathcal{HD}}((u\,\la\,a).\underline v)\ .
\end{eqnarray}
First we have, using \textit{(ii)} for $n$,
\begin{eqnarray} \label{cvbyudiskify}
\bar\partial_{\mathcal{HD}}( u\bullet(a.\underline v)) &=& ((\pi^{0,1}\tens\id)\square\,u)\bullet(a. \underline v)
- (\ev\tens \id)(u\tens \phi(\pi^{0,1}\wedge\pi^{1,0})\square(a.\xi))\tens \underline w \cr
&&  - \ (\ev\tens \id)(u\tens \phi(a.\xi\wedge\eta_i))\tens \eta^i\bullet \underline w \cr
&& -\  (\ev\tens \id)(u\tens \phi(\pi^{0,1}\wedge\pi^{1,0})\square(\bar\partial a))\tens \underline v \cr
&&  - \ (\ev\tens \id)(u\tens \phi(\bar\partial a\wedge\eta_i))\tens \eta^i\bullet \underline v \ ,\cr
\bar\partial_{\mathcal{HD}}( (u.a)\bullet\underline v) &=& ((\pi^{0,1}\tens\id)\square(u.a))\bullet \underline v
- (\ev\tens \id)(u.a\tens \phi(\pi^{0,1}\wedge\pi^{1,0})\square(\xi))\tens \underline w \cr
&&  - \ (\ev\tens \id)(u.a\tens \phi(\xi\wedge\eta_i))\tens \eta^i\bullet \underline w\ ,\cr
 \bar\partial_{\mathcal{HD}}((u\,\la\,a).\underline v) &=& \bar\partial(u\,\la\,a)\tens \underline v+(u\,\la\,a).\xi\tens\underline w \cr
 &=&  (((\pi^{0,1}\tens\id)\square\, u)\,\la\,a)\tens \underline v + (\ev\tens\id)(u\tens(\id\tens \pi^{0,1})\square\partial a) \tens \underline v \cr && +\ (u\,\la\,a).\xi\tens\underline w  \ .
\end{eqnarray}
For the last line of (\ref{cvbyudiskify}) we used
\begin{eqnarray}
\bar\partial(u\,\la\,a) &=& \bar\partial\,\ev(u\tens\partial a)\cr
&=& ((\pi^{0,1}\tens\id)\square\, u)\,\la\,a + (\ev\tens\id)(u\tens(\id\tens \pi^{0,1})\square\partial a)\ .
\end{eqnarray}
Now from (\ref{cvbyudiskify}),
\begin{eqnarray} \label{cvbysy}
&&\bar\partial_{\mathcal{HD}}( u\bullet(a.\underline v)) -\bar\partial_{\mathcal{HD}}( (u.a)\bullet\underline v) 
-  \bar\partial_{\mathcal{HD}}((u\,\la\,a).\underline v)   \cr 
&=& ((\pi^{0,1}\tens\id)\square\,u)\bullet(a. \underline v) - ((\pi^{0,1}\tens\id)\square(u.a))\bullet \underline v 
- (((\pi^{0,1}\tens\id)\square\, u)\,\la\,a)\tens \underline v \cr
&& - \  (\ev\tens \id)(u\tens \phi(\pi^{0,1}\wedge\pi^{1,0})\square(a.\xi))\tens \underline w +(\ev\tens \id)(u.a\tens \phi(\pi^{0,1}\wedge\pi^{1,0})\square(\xi)).\underline w \cr
&& -\  (\ev\tens \id)(u\tens \phi(\pi^{0,1}\wedge\pi^{1,0})\square(\bar\partial a))\tens \underline v 
- (\ev\tens\id)(u\tens(\id\tens \pi^{0,1})\square\partial a) \tens \underline v\cr
&&  - \ (\ev\tens \id)(u\tens \phi(\bar\partial a\wedge\eta_i))\tens \eta^i\bullet \underline v - (u\,\la\,a).\xi\tens\underline w  \cr
&=& - \ ((\pi^{0,1}\tens\pi^{1,0})\sigma(u\tens\extd a))\bullet \underline v   -  (\ev\tens \id)(u\tens \phi(\bar\partial a\wedge\eta_i))\tens \eta^i\bullet \underline v    \cr
&& - \  (\ev\tens \id)(u\tens \phi(\pi^{0,1}\wedge\pi^{1,0})\sigma^{-1}(\extd a\tens\xi))\tens \underline w  - (u\,\la\,a).\xi\tens\underline w  \cr
&& -\  (\ev\tens \id)(u\tens \phi(\pi^{1,1}(\wedge\square)(\bar\partial a))\tens \underline v 
- (\ev\tens\id)(u\tens\phi(\pi^{1,1}(\wedge\square)\partial a) \tens \underline v \cr
&=& - \big(((\pi^{0,1}\tens\pi^{1,0})\sigma(u\tens\extd a))   +  (\ev\tens \id)(u\tens \phi(\bar\partial a\wedge\eta_i))\tens \eta^i\big)\bullet \underline v    \cr
&& - \  \big((\ev\tens \id)(u\tens \phi(\pi^{0,1}\wedge\pi^{1,0})\sigma^{-1}(\extd a\tens\xi))  + (u\,\la\,a).\xi\big)\tens\underline w  \cr
&& -\  (\ev\tens \id)(u\tens \phi(\pi^{1,1}(\wedge\square)(\extd a))\tens \underline v   \ .
\end{eqnarray}
The last line of (\ref{cvbysy}) vanishes by our assumption on the torsion. The expression in (\ref{cvbysy}) vanishing would be implied by the following equations (\ref{kjasyfgvkuy}):
\begin{eqnarray}  \label{kjasyfgvkuy}
\partial a\tens\xi &=& -\  \phi(\pi^{0,1}\wedge\pi^{1,0})\sigma^{-1}(\extd a\tens\xi)\ ,\cr
(\pi^{0,1}\tens\pi^{1,0})\sigma(u\tens\extd a)   &=& -\   (\ev\tens \id)(u\tens \phi(\bar\partial a\wedge\eta_i))\tens \eta^i\ .
\end{eqnarray}
Remember that $\xi\in\Omega^{0,1}A$ in (\ref{kjasyfgvkuy}). Using the formula (\ref{sigrefdefdual}) linking $\sigma$ and $\sigma^{-1}$, the last equation in (\ref{kjasyfgvkuy}) is equivalent to
\begin{eqnarray} \label{isdycguys}
(\pi^{1,0}\tens\pi^{0,1})\sigma^{-1}(\extd a\tens \kappa) \ =\ -\,\phi(\bar\partial a\wedge \kappa)\ ,
\end{eqnarray}
for $\kappa\in\Omega^{1,0}A$. As we have assumed that $\wedge:\Omega^{1,0}A\tens_A \Omega^{0,1}A\to \Omega^{1,1}A$ is an isomorphism with inverse $\phi$, we can apply $\wedge$ to the first equation in (\ref{kjasyfgvkuy}) and to (\ref{isdycguys}) to get two equivalent equations,
\begin{eqnarray} \label{vcisuyiu}
\partial a\wedge\xi &=& - \  (\pi^{0,1}\wedge\pi^{1,0})\sigma^{-1}(\extd a\tens\xi)\ =\ -\ \pi^{1,1}\wedge\,\sigma^{-1}(\extd a\tens\xi)\cr
\bar\partial a\wedge \kappa &=& -\ (\pi^{1,0}\wedge\pi^{0,1})\sigma^{-1}(\extd a\tens \kappa)\ =\ -\ \pi^{1,1}\wedge\,\sigma^{-1}(\extd a\tens \kappa)\ ,
\end{eqnarray}
and this is given by \textit{(d)} in the hypothesis, which finishes verifying well definition. 

Now verifying \textit{(ii)} for $n+1$ is rather easier:
\begin{eqnarray}
\bar\partial_{\mathcal{HD}}( (a.u)\bullet\underline v) &=& ((\pi^{0,1}\tens\id)\square(a.u))\bullet \underline v
- (\ev\tens \id)(a.u\tens \phi(\pi^{0,1}\wedge\pi^{1,0})\square(\xi))\tens \underline w \cr
&&  - \ (\ev\tens \id)(a.u\tens \phi(\xi\wedge\eta_i))\tens \eta^i\bullet \underline w \cr
&=& \bar\partial a\tens u\bullet  \underline v +a.  \bar\partial_{\mathcal{HD}}( u\bullet\underline v) \ .\quad\blacksquare
\end{eqnarray}

\begin{propos} \label{cbhadjsvhjgcgfxvdfsz}
Suppose that the conditions of Proposition~\ref{cbhadjsvhjgcgfx} are satisfied. If $\underline v\in \mathcal{T}\mathrm{Vec}^{*,0}A$ and $a\in A$ with $\bar\partial a=0$, then
$\bar\partial(\underline v\,\la\,a)=\bar\partial_{\mathcal{HD}}(\underline v)\,\la\,a$. 
\end{propos}
\noindent {\bf Proof:}\quad The statement is true for $n=1$ by Lemma~\ref{kuyfdtyytr44s}. Now assume that it us true for $n$. 
To show the statement for $n+1$, we set $\bar \partial a=0$, and then, setting $b=\underline v\,\la\,a$ and using Lemma \ref{kuyfdtyytr44s},
\begin{eqnarray*} 
\bar\partial((u\bullet\underline v)\,\la\,a) &=& \bar\partial(u\,\la\,b) \cr
&=& (\pi^{0,1}\tens\ev)(\square\,u\tens \partial b)
- (\ev\tens \id)(u\tens \phi\,\pi^{1,1}\wedge\square\,\bar\partial b)\ .
\end{eqnarray*}
Using the statement for $n$, and setting $\xi\tens \underline w=\bar\partial_{\mathcal{HD}}( \underline v)$, we have $\bar\partial b = \xi.(\underline w\,\la\,a)$, so
\begin{eqnarray*}
(\id\tens\pi^{1,0})\square\,\bar\partial b &=& (\id\tens\pi^{1,0})\square(\xi).(\underline w\,\la\,a)+
\xi\tens\partial(\underline w\,\la\,a)\ ,
\end{eqnarray*}
and substituting this gives 
\begin{eqnarray} \label{cbksuyjkhg}
\bar\partial((u\bullet\underline v)\,\la\,a)
&=& (\pi^{0,1}\tens\ev)(\square\,u\tens \partial b)
- (\ev\tens \id)(u\tens \phi(\pi^{0,1}\wedge\pi^{1,0})\square(\xi)).(\underline w\,\la\,a)   \cr
&& - \ (\ev\tens \id)(u\tens \phi\,\pi^{1,1}(\xi\wedge\partial(\underline w\,\la\,a)))\ .
\end{eqnarray}
The first two terms of (\ref{cbksuyjkhg}) correspond to the first two terms of the displayed equation in the statement of 
Proposition~\ref{cbhadjsvhjgcgfx}. For the last term, we use a dual basis $\eta_i\tens\eta^i\in\Omega^{1,0}A\tens_A \mathrm{Vec}^{1,0}A$, and then
\begin{eqnarray*}
(\ev\tens \id)(u\tens \phi\,\pi^{1,1}(\xi\wedge\partial(\underline w\,\la\,a))) &=& 
(\ev\tens \id)(u\tens \phi\,\pi^{1,1}(\xi\wedge\eta_i)).\eta^i(\partial(\underline w\,\la\,a))\ .
\end{eqnarray*}
giving the last line in Proposition~\ref{cbhadjsvhjgcgfx}. \quad$\blacksquare$

\medskip The rather unpleasant expression for $\bar\partial_{\mathcal{HD}}$ in Proposition~\ref {cbhadjsvhjgcgfx} will become somewhat clearer if we refer to
Examples~\ref{iouycgziuyf1} and \ref{iouycgziuyf2}. In fact we can rewrite the recursive definition as
\begin{eqnarray}  \label{vyuiutyc}
\bar\partial_{\mathcal{HD}}( u\bullet\underline v) &=& \bar\partial u\bullet \underline v
+u\,\la\,(\bar\partial_{\mathcal{HD}}(\underline v))\ .
\end{eqnarray}
The first term in (\ref{vyuiutyc}) obviously is the first term in the recursive definition in the statement of Proposition~\ref {cbhadjsvhjgcgfx}, but seeing that the second term of  (\ref{vyuiutyc})  is the second and third terms in the recursive definition is a little more difficult. In Examples~\ref{iouycgziuyf1} and \ref{iouycgziuyf2} we see that we have a left $\partial$ covariant derivative on $\Omega^{0,1}A\tens_A \mathcal{T}\mathrm{Vec}^{*,0}A_\bullet$, where we use the fact that we actually have a bimodule covariant derivative on $\Omega^{0,1}A$, and this gives the left action of $u$. To get the equations to match, we have to remember that $\pi^{1,1}$ composed with the torsion of $\square$ vanishes, as given in the conditions in Proposition~\ref {cbhadjsvhjgcgfx}. 

\begin{propos} \label{gvuyiiuyf}
For all $\underline v,\underline w\in \mathcal{T}\mathrm{Vec}^{*,0}A_\bullet$, 
\begin{eqnarray*} 
\bar\partial_{\mathcal{HD}}( \underline v\bullet\underline w) &=& \bar\partial_{\mathcal{HD}}(\underline v)\bullet \underline w
+\underline v\,\la\,(\bar\partial_{\mathcal{HD}}(\underline w))\ .
\end{eqnarray*}
Thus the $\bullet$ product of holomorphic differential operators is holomorphic. 
\end{propos}
\noindent {\bf Proof:}\quad The proof is by induction on $n$, where $\underline v\in \mathrm{Vec}^{\tens n,0}A$. The 
$n=1$ case is given by (\ref{vyuiutyc}). Suppose the statement works for all values $\le n$. Now, for 
$u\in \mathrm{Vec}A$ and $\underline v\in \mathrm{Vec}^{\tens n,0}A$,
\begin{eqnarray*} 
\bar\partial_{\mathcal{HD}}( (u\bullet \underline v)\bullet\underline w) &=&
\bar\partial_{\mathcal{HD}}( u\bullet (\underline v\bullet\underline w))\cr
&=& 
 \bar\partial_{\mathcal{HD}}(u)\bullet(\underline v\bullet\underline w)
+u\,\la\,(\bar\partial_{\mathcal{HD}}(\underline v\bullet\underline w))\cr
&=& 
 \bar\partial_{\mathcal{HD}}(u)\bullet(\underline v\bullet\underline w)
+u\,\la\,(\bar\partial_{\mathcal{HD}}(\underline v)\bullet \underline w)
+u\,\la\,(\underline v\,\la\,(\bar\partial_{\mathcal{HD}}(\underline w)))   \cr
&=& 
( \bar\partial_{\mathcal{HD}}(u)\bullet\underline v)\bullet\underline w
+(u\,\la\,\bar\partial_{\mathcal{HD}}(\underline v))\bullet \underline w
+(u\bullet\underline v)\,\la\,(\bar\partial_{\mathcal{HD}}(\underline w))   \ .
\end{eqnarray*}
To get to the last line of this equation we have used, in order, the associativity of $\bullet$, the fact that the left $\partial$ covariant derivative in Example~\ref{iouycgziuyf1} is a right $\mathcal{T}\mathrm{Vec}^{*,0}A_\bullet$ module map, and the fact that $\la$ is an action of the algebra $\mathcal{T}\mathrm{Vec}^{*,0}A_\bullet$. \quad$\blacksquare$

\begin{defin}
Suppose that the conditions of Proposition~\ref{cbhadjsvhjgcgfx} are satisfied.
An element $\underline v\in \mathcal{T}\mathrm{Vec}^{*,0}A_\bullet$ is called holomorphic if
$\bar\partial_{\mathcal{HD}}(\underline v)=0$. 
From Proposition~\ref{cbhadjsvhjgcgfxvdfsz} a holomorphic $\underline v\in \mathcal{T}\mathrm{Vec}^{*,0}A_\bullet$ acting on a holomorphic $a\in A$ gives a holomorphic element of $A$. 
From Proposition~\ref{gvuyiiuyf} the holomorphic elements are closed under the $\bullet$ product. 
\end{defin}

\end{document}